\documentclass[a4paper]{amsart}

\usepackage{amssymb}
\usepackage{amsthm}
\usepackage{amsmath} \usepackage{latexsym}
\usepackage{color}%

\newtheorem{thm}{Theorem}[section]
\newtheorem{prop}[thm]{Proposition}
\newtheorem{lemma}[thm]{Lemma}
\newtheorem{cor}[thm]{Corollary}
\theoremstyle{definition}

\newtheorem{ex}[thm]{Example}
\newtheorem{exs}[thm]{Examples} \newtheorem{rem}[thm]{Remark}

 \DeclareMathOperator{\Spec}{Spec}
 \DeclareMathOperator{\Max}{Max}
 
 \DeclareMathOperator{\Na}{Na}
 \DeclareMathOperator{\Kr}{Kr}

 \DeclareMathOperator{\Sup}{Sup}
  \DeclareMathOperator{\di}{dim}

  \DeclareMathOperator{\htt}{ht}
\begin{document}

\title[SEMISTAR DEDEKIND DOMAINS]{SEMISTAR DEDEKIND DOMAINS}
\author[ ]{Said El Baghdadi$^1$, Marco Fontana$^2$, Giampaolo Picozza$^2$}

 \dedicatory{$^1$ Department of Mathematics, Facult\'e des Sciences et Techniques,\\
P.O. Box 523, Beni Mellal, Morocco \\ \rm \; \texttt{baghdadi@fstbm.ac.ma} \rm \\
\vskip 0.3cm \textit{$^2$Dipartimento di Matematica, Universit\`a
degli
Studi ``Roma Tre'', \\
 Largo San Leonardo Murialdo 1, 00146 Roma} \\ \rm  \texttt {fontana@mat.uniroma3.it \;
 picozza@mat.uniroma3.it}}

 \thanks{\it 2000 Mathematics Subject Classification: \rm 13A15, 13G05, 13E99.}%
 \thanks{\it Acknowledgment: \rm Part of this work was done while the 
 first named author was visiting the  Mathematical Department of 
 Universit\`a degli Studi ``Roma Tre'', with a visiting grant by INdAM.  The second and the third named authors were supported in 
part by a research grant by MIUR 2003/2004.}%

\keywords{Dedekind  domain, semistar operation, Krull domain, Mori 
domain, Pr\"ufer domain, Nagata ring.}%


\begin{abstract}  
Let $D$ be an integral domain and $\star$ a
semistar operation on $D$. As a generalization of the notion of Noetherian
domains to the semistar setting,  we say that $D$ is a
$\star$--Noetherian domain if it has the ascending chain condition
on the set of its quasi--$\star$--ideals. 
On the other hand, as an
extension the notion of Pr\"ufer domain (and of Pr\"{u}fer
$v$--multiplication domain),
we say that $D$ is a Pr\"ufer $\star$--multiplication domain (P$\star$MD, for short) if $D_M$ is a
valuation domain,  for each quasi--$\star_{_{\!f}}$--maximal ideal $M$ of $D$. Finally,
recalling that a Dedekind domain is a Noetherian
Pr\"{u}fer domain,
we define a $\star$--Dedekind domain to be an integral domain
which is $\star$--Noetherian and a P$\star$MD.  For the identity
semistar operation $d$, this definition coincides with that of the
usual Dedekind domains and when the semistar operation is the
$v$--operation, this notion gives rise to Krull domains. Moreover,
Mori domains not strongly Mori are $\star$--Dedekind for a suitable 
spectral semistar operation.\\
Examples show that $\star$--Dedekind domains are not necessarily
integrally closed nor one-dimensional, although they mimic
various aspects, varying  according to  the choice of $\star$, of
the ``classical'' Dedekind domains.
In any case, a $\star$--Dedekind domain is an
integral domain $D$ having a   Krull overring $T$ (canonically
associated to $D$ and $\star$) such that
the semistar operation $\star$ is essentially ``univocally associated'' to the
$v$--operation on $T$.

In the present paper, after a preliminary study of  $\star$--Noetherian 
domains, we investigate 
the 
$\star$--Dedekind
domains.
We extend to the $\star$--Dedekind domains the main
classical results and several characterizations   proven   for
Dedekind domains. In particular,  we obtain a characterization of
a $\star$--Dedekind domain by a property of decomposition of any 
semistar ideal into a
``semistar product'' of prime ideals.
Moreover, we show that  an integral domain $D$ is a $\star$--Dedekind domain if
and only if the Nagata semistar domain $\Na(D, \star)$ is a Dedekind
domain. Several applications of the general results are given for 
special cases of the semistar operation $\star$.  
\end{abstract}

\maketitle

\date{\today}

\section{Introduction and background results}

Dedekind domains play a crucial role in classical algebraic number
theory and their study gave a relevant contribution  to a rapid development
of commutative ring theory and ideal theory: Noetherian, Krull and Pr\"{u}fer domains
arose in the early stages of these theories,
for generalizing different aspects of Dedekind domains.
%


Star operations provided new insight in multiplicative ideal theory.
For instance, the
use of the  $v$-- and $t$--operations has
produced a common
treatment
and a deeper understanding of Dedekind and Krull domains.
In 1994, Okabe
and  Matsuda  \cite{Okabe/Matsuda:1994} introduced the semistar operations, extending
the notion of star operation and the related classical theory
of ideal systems, based on the pioneering works by W. Krull, E. Noether, H.
Pr\"ufer and P. Lorenzen (cf. \cite{Jaffard} and \cite{HK98}).
Semistar operations,   due to a major grade
of flexibility with respect to star operations,  provide a natural and
general setting for a wide class  of questions and  for  a deeper and comparative study of
several relevant classes of integral domains (cf. for instance
\cite{Okabe/Matsuda:1994}, \cite{Matsuda/Sato:1996},
\cite{Matsuda/Sugatani:1995}, \cite{FH2000}, \cite{FL01}, \cite{FL01b}, \cite{FL03},
\cite{FJS03}, and \cite{Koch:2000}).

In this paper, we explore a general theory of Dedekind--type domains,
depending on a semistar operation. A first attempt in this direction
was done by Aubert \cite{Au}
and, more extensively, by Halter-Koch \cite[Chapter 23]{HK98},
where the author investigated Dedekind
domains in the language of finitary ideal systems on commutative
monoids. Our approach is
based on the classical
multiplicative ideal theory on integral domains,
as in Gilmer's book \cite{Gil92},
extended in a natural way to the semistar case. This approach has already produced
a generalized and covenient setting for
considering Kronecker function rings (\cite{FL01}, \cite{FL01b}, and
\cite{FL03}), Nagata
rings \cite{FL03}, and Pr\"ufer multiplication domains \cite{FJS03}.\\
Note that the module
systems approach on commutative monoids, developed recently by Halter-Koch
in \cite{Koch:JA}, provides  an alternative general frame for
(re)con\-si\-de\-ring semistar
operations on integral domains and related topics.  More precisely, 
most of the results contained in this paper are of purely 
multiplicative nature and remain valid in the more general setting of 
commutative cancellative monoids (cf. also Remark \ref{new:rk}).  

\medskip
Recall that a Dedekind domain is a Noetherian Pr\"{u}fer domain. Let $D$ be an integral domain and $\star$ a semistar operation on
$D$. As a generalization of Noetherian domains to the semistar
setting, we define $D$ to be a $\star$--Noetherian domain if it has
the ascending chain condition on the set of the ideals of $D$
canonically associated to $\star$ (called quasi--$\star$--ideals); equivalently,
a $\star$--Noetherian domain is a domain in which each nonzero ideal is
$\star_{_{\!f}}$--finite (Lemma \ref{lemma:3.1} and Remark \ref{rk:3.6} (1)).
For instance, as we will see later, Noetherian,\ Mori, and strong
Mori domains are examples of $\star$--Noetherian domains, for
different $\star$--operations.

On the other hand, as an extension the notion of Pr\"ufer domain (and 
of 
Pr\"{u}fer $v$--multiplication domain), given a semistar operation
$\star$  on an integral domain
$D$,  we say that $D$ is a Pr\"ufer
semistar multiplication domain (P$\star$MD, for short) if $D_M$ is a
valuation domain,  for each maximal element $M$ in the nonempty set of the ideals
of $D$ associated to the finite type semistar operation canonically
deduced from $\star$ (i.e.,  the quasi--$\star_{_{\!f}}$--maximal ideal of $D$).
Finally, we
define a $\star$--Dedekind domain ($\star$--DD, for short) to be an integral domain which is $\star$--Noetherian and
a P$\star$MD.  For the identity semistar operation $d$, this
definition coincides with that of the usual Dedekind domains and when
the semistar operation is the $v$--operation, this notion gives rise to Krull
domains.  Moreover,
Mori domains not strongly Mori are $\star$--Dedekind for a suitable 
spectral semistar operation  (Example \ref{ex:4.22}).

In the general semistar setting, $\star$--Dedekind domains are not necessarily
integrally closed nor one-dimensional, although they mimic
various aspects, varying  according to  the choice of $\star$, of
the ``classical'' Dedekind domains.
In any case, a $\star$--Dedekind domain is an
integral domain $D$ having a   Krull overring $T$ (canonically
associated to $D$ and $\star$) such that
the semistar operation $\star$ is essentially ``univocally associated'' to the
$v$--operation on $T$ (Remark \ref{rk:4.17}). 
 
In the present paper we develop a theory which
enlightens different facets of the $\star$--Dedekind domains and shows
the interest in studying these new classes of integral domains of
Dedekind-type, parametrized by semistar operations.  After recalling
in the present section the
main data needed for this work,  in Section 2, as a first step to the main goal,
we introduce and study the concept of  ``semistar almost
Dedekind domains'' ($\star$--ADD, for short), which provides a natural generalization of the classical
notion of almost Dedekind domains.  Our study, in the particular case
of $\star = v$,  extends and
completes the investigation on $t$--almost Dedekind domains initiated
by Kang \cite[Section IV]{Ka89}.  Among the main results  proven  
 in this
section, we have that an integral domain $D$ is a $\star$--ADD if and only
if the Nagata semistar domain $\Na(D, \star)$ is an almost
Dedekind domain (in particular,  in this case, $\Na(D, \star)$
coincides with the Kronecker semistar function ring $\Kr(D, \star)$).

Section 3
is devoted to the study of the semistar Noetherian domains. In
particular, we investigate the local--global
behaviour of this notion and we obtain several relevant results on
$\star$--Noetherian domains, in case of
stable semistar operations.

In Section 4, we introduce and study the semistar Dedekind
domains. We extend to the $\star$--Dedekind domains the main classical
results and several characterizations   proven   for Dedekind domains. In
particular,  we obtain a
characterization of  a $\star$--Dedekind domain by a property of
decomposition of any semistar ideal into a ``semistar product'' of prime ideals.
Moreover, we
show that  an integral domain $D$ is a $\star$--DD if and only
if the Nagata semistar domain $\Na(D, \star)$ is a
Dedekind domain (in particular,  in this case, $\Na(D, \star)$
coincides with the Kronecker semistar function ring $\Kr(D, \star)$, 
 which    is in fact a PID).

\bigskip

\centerline{ $\ast \quad \ast \quad  \ast \quad \ast \quad  \ast $}
\bigskip

Let $D$ be an integral domain with quotient field $K$. Let
$\boldsymbol{\overline{F}}(D)$
   denote    the set of all nonzero $D$--submodules of $K$  and let
  $\boldsymbol{F}(D)$ represent the nonzero fractional ideals of $D$ (i.e.,
  $\boldsymbol{F}(D):=\{E \in \boldsymbol{\overline{F}}(D)\mid dE \subseteq
  D \mbox{ for some nonzero element } d \in D\} $).  Finally, let
  $\boldsymbol{f}(D)$  be the set of all nonzero finitely
  generated $D$-submodules of $K$ (it is clear that $\boldsymbol{f}(D)
  \subseteq \boldsymbol{F}(D)    \subseteq
  \boldsymbol{\overline{F}}(D)$).

A map $\star: \boldsymbol{\overline{F}}(D) \to
\boldsymbol{\overline{F}}(D), E \mapsto E^\star$, is called a
\emph{semistar operation} on $D$ if, for all $x \in K$, $x \neq
0$, and for all $E,F \in \boldsymbol{\overline{F}}(D)$, the
following properties hold:
\begin{enumerate}
\item[$(\star_1)$] $(xE)^\star=xE^\star$;
\item[$(\star_2)$] $E \subseteq F$ implies $E^\star \subseteq
F^\star$;
\item[$(\star_3)$] $E \subseteq E^\star$ and $E^{\star \star} :=
\left(E^\star \right)^\star=E^\star$;
\end{enumerate}
cf.  for instance \cite{Okabe/Matsuda:1994}, \cite{Matsuda/Sugatani:1995},
\cite{Matsuda/Sato:1996} and \cite{FH2000}.  Recall that,  given a
semistar operation $\star$ on $D$, for all $E,F \in
\boldsymbol{\overline{F}}(D)$,  the following basic formulas follow
easily from the axioms:    $$\begin{array}{rl} (EF)^\star =& \hskip -7pt
(E^\star F)^\star =\left(EF^\star\right)^\star =\left(E^\star
F^\star\right)^\star\,;\\
(E+F)^\star =& \hskip -7pt \left(E^\star + F\right)^\star= \left(E
+
F^\star\right)^\star= \left(E^\star + F^\star\right)^\star\,;\\
(E:F)^\star \subseteq & \hskip -7pt (E^\star :F^\star) = (E^\star
:F) =
\left(E^\star :F\right)^\star,\;\, \mbox{\rm if \ } (E:F) \neq 0;\\
(E\cap F)^\star \subseteq & \hskip -7pt E^\star \cap F^\star =
\left(E^\star \cap F^\star \right)^\star,\;\, \mbox{\rm if \ } E\cap F \neq
(0)\,;
\end{array}
$$
 \noindent cf.  for instance \cite[Theorem 1.2 and p.  174]{FH2000}.

A \emph{(semi)star operation} $\star$ on an integral domain $D$ is
a semistar operation, that restricted to the set
$\boldsymbol{F}(D)$ of fractional ideals, is a  \emph{star operation} on $D$
\cite[(32.1)]{Gil92}.     It is very easy to see that  a semistar
operation $\star$ on $D$    is a (semi)star operation  (on $D$)
if and only if $D^\star=D$.

\begin{exs} \label{ex:1.1} \bf (1) \rm The first (trivial) examples of
semistar operations are given by $d_{D}$ (or, simply, $d$), called the
\it identity (semi)star operation on \rm $D$, defined by $E^{d_{D}} :=
E$, for each $E \in \boldsymbol{\overline{F}}(D)$ and by
$e_{D}$   (or,
simply,   $e$),    defined by   $E^{e_{D}}    := K$,    for each $E \in
\boldsymbol{\overline{F}}(D)$.\\
More generally, if $T$ is an overring of
$D$, we can define a semistar operation on $D$, denoted by
$\star_{\{T\}}$ and defined by $E^{\star_{\{T\}}} := ET$, for each $E
\in \boldsymbol{\overline{F}}(D)$.  It is obvious that $d_{D} =
\star_{\{D\}}$,   $e_{D}    = \star_{\{K\}}$    and that $\star_{\{T\}}$ is a
semistar non--(semi)star operation on $D$ if and only if $D\subsetneq
T$.

\bf (2) \rm If $\star$ is a semistar operation on $D$, then we can
consider a map\ $\star_{\!_f}: \boldsymbol{\overline{F}}(D) \to
\boldsymbol{\overline{F}}(D)$ defined, for each $E \in
\boldsymbol{\overline{F}}(D)$, as follows:

\centerline{$E^{\star_{\!_f}}:=\bigcup \{F^\star\mid \ F \in \boldsymbol{f}(D)
\mbox{ and } F \subseteq E\}$.}

\noindent It is easy to see that $\star_{\!_f}$ is a semistar operation on $D$,
called \emph{the semistar operation of finite type associated to
$\star$}.  Note that, for each $F \in \boldsymbol{f}(D)$,
$F^\star=F^{\star_{\!_f}}$.  A semistar operation $\star$ is called a
\emph{semistar operation of finite type} if $\star=\star_{\!_f}$.  It is
easy to see that $(\star_{\!_f}\!)_{\!_f}=\star_{\!_f}$ (that is,
$\star_{\!_f}$ is of finite type). \\
 For instance, if $v_{D}$ (or, simply, $\,v\,$) is \it the
$v$--(semi)star operation on $D$ \rm defined by $E^v := (E^{-1})^{-1},$
for each $E \in \overline{\boldsymbol{F}}(D)\,,$ with $E^{-1} := (D
:_{\mbox{\tiny \it K}} E) := \{ z \in K \; | \;\; zE \subseteq D
\}\,$,\, then the
semistar operation of finite type $(v_{D})_{f}$ (or, simply,
$v_{_{\!f}}$) associated to $v_{D}$ is denoted by $t_{D}$ (or, simply,
$t$) and it is called \it the $t$--(semi)star operation on $D$\, \rm
(note that $D^v =D^t = D$).  \\
Note also that, for each overring $T$ of $D$, the semistar operation
$\star_{\{T\}}$ on $D$ is a semistar operation of finite type.

\bf (3) \rm If $\Delta \subseteq \Spec(D)$, the map $\star_{\Delta}:
\boldsymbol{\overline{F}}(D)  \rightarrow  \boldsymbol{\overline{F}}(D)$, $E
\mapsto E^{\star_\Delta} := \bigcap \{ED_P \mid \  P \in \Delta \}$, is a
semistar operation on $D$ \cite[Lemma 4.1]{FH2000}.    A semistar operation
$\star$ is called  a \emph{spectral semistar operation on} $D$    if
there exists a subset $\Delta$ of $\Spec(D)$ such that $\star =
\star_\Delta$.  If $\Delta=\{P\}$, then $\star_{\{P\}}$ is the spectral
semistar operation on $D$ defined by $E^{\star_{\{P\}}}:=ED_P$, for each
$E \in \boldsymbol{\overline{F}}(D)$,  i.e. $\star_{\{P\}}
=\star_{\{D_{P}\}}$.     If  $\Delta= \emptyset$, then
$\star_{\emptyset} = e_{D}$.  \\
We say that a semistar operation is
\emph{stable (with respect to   finite    intersections)} if $(E \cap F)^\star=
E^\star \cap F^\star$, for each $E,F \in \boldsymbol{\overline{F}}(D)$.
For a spectral semistar operation the following properties hold \cite[Lemma
4.1]{FH2000}:
\begin{enumerate}
   \item[\bf{(3.a)}] For each $E\in \boldsymbol{\overline{F}}(D)$ and for
   each $P\in \Delta\,$, $ED_P=E^{\star_\Delta}D_P$.
   \item[\bf{(3.b)}] The semistar operation $\star_\Delta$ is stable.
   \item[\bf{(3.c)}] For each
   $P\in \Delta$, $P^{\star_\Delta}\cap D = P$.
   \item[\bf{(3.d)}] For each nonzero integral ideal $I$ of $D$ such that
   $I^{\star_\Delta}\cap D\neq D$, there exists a prime ideal $P\in \Delta$
   such that $I\subseteq P$.  \end{enumerate}
\end{exs}

If ${\star_1}$ and ${\star_2}$ are two semistar operations on $D$,
we  set\   ${\star_1} \leq {\star_2}$,\ if $E^{\star_1} \subseteq
E^{\star_2}$, for each $E \in \boldsymbol{\overline{F}}(D)$.  It is easy
to see that $\star_1 \leq \star_2$ if and only if
$\left(E^{\star_{1}}\right)^{\star_{2}} = E^{\star_2}=
\left(E^{\star_{2}}\right)^{\star_{1}}$.   Obviously, if ${\star_1}
\leq {\star_2}$, then ${(\star_1)}_{\!_f} \leq {(\star_2)}_{\!_f} $;
moreover, for each semistar operation $\star$  on $D$,    we have
 $d_{D} \leq \star_{\!_f} \leq \star \leq  e_{D}$.   In particular,
$t_{D}\leq v_{D}$; furthermore, it is not difficult to see that, for
each (semi)star operation $\star$ on $D$, we have $\star \leq v_{D}$ and
$\star_{_{\!f}}\leq t_{D}$ \cite[Theorem 34.1(4)]{Gil92}.

 A \it quasi--$\star$--ideal \rm of $D$ is a  nonzero     ideal $I$ of
 $D$ such that $I = I^\star \cap D$.  This notion generalizes the notion of
 \it $\star$--ideal \rm for a star operation on $D$,  which is a nonzero
 ideal $I$ of $D$ such that $I = I^\star$.     More precisely, it is clear
 that, for a (semi)star operation  $\star$, the quasi--$\star$--ideals
 coincide with the $\star$--ideals.  \\
 Note that each  nonzero    ideal $I$ of $D$, such that $I^\star \subsetneq
 D^\star$, is contained in a (non trivial) quasi--$\star$--ideal  of
 $D$:    in fact, the ideal $I^\star \cap D$ is a quasi--$\star$--ideal
 of $D$ and $I \subseteq I^\star \cap D$.

A \emph{quasi--$\star$--prime}  of $D$    is a  nonzero    prime
ideal  of $D$    that is also a quasi--$\star$--ideal  of $D$.     A
\it quasi--$\star$--maximal ideal \rm  of $D$    is a (proper) ideal
 of $D$, which    is maximal in the set of all quasi--$\star$--ideals
of $D$.\\
If $\star$ is a semistar operation of finite type  on $D$,  with
$D \neq K$,    each
quasi--$\star$--ideal of $D$ is contained in a quasi--$\star$--maximal
ideal.  Moreover, each quasi--$\star$--maximal ideal of $D$ is prime
\cite[Lemma 4.20]{FH2000}.  We denote by $\mathcal{M}(\star)$ the set of
all quasi--$\star$--maximal ideals of $D$.   Thus, if $\star =
\star_{_{\!f}}$ and $D$ is not a field, then $\mathcal{M}(\star) \neq
\emptyset$.

\setcounter{thm}{0}
\begin{exs} \label{ex:1,4-6} \bf (4) \rm If $\star$ is a semistar
operation on $D$, we denote by ${\tilde{\star}}$ the spectral semistar
operation $\star_{\mathcal{M}(\star_{\!_f})}$, induced by the set
$\mathcal{M}(\star_{\!_f})$ of the quasi--$\star_{\!_f}$--maximal ideals
of $D$, i.e. for each $E\in \boldsymbol{\overline{F}}(D)$:

\centerline{$E^{\tilde{\star}} := \bigcap \{ED_{Q }\mid \ Q\in
\mathcal{M}(\star_{\!_f})\} $\,.   }

\noindent The semistar operation ${\tilde{\star}}$ is stable and of finite type
and ${\tilde{\star}} \leq \star_{\!_f}$ (cf.   \cite{FL01} and
also   \cite[p.  181]{FH2000}
for an equivalent definition of ${\tilde{\star}}$;   see   \cite{Hedstrom/Houston:1980},
\cite{Glaz/Vasconcelos:1977},  \cite[Section
2]{AC2000} for an analogous construction in the star setting).   Note
that, when $\star$ is a (semi)star operation on $D$, then also
$\tilde{\star}$ is a (semi)star operation on $D$.     \\
 If $\star = d_{D}$, then obviously ${\tilde{\star}} =d_{D}$.  If $\star = v_{D}$,  then ${\tilde{\star}}$ is the (semi)star operation
that we denote  by $w_{D}$ (or, simply, $w$), following    Wang Fanggui
and R.L. McCasland (cf.  \cite{WMc97},
\cite{WMc99} and \cite{Fanggui:2001}).  Note that, for $\Delta=
\mathcal{M}(\star_{_{\!f}})$, the semistar operation $\tilde{\star}$
satisfies the properties {\bf (3.a)}--{\bf (3.d)}\rm, stated above for a
general spectral semistar operation.

\bf (5) \rm Let $D$ be an integral domain and $T$ an overring of $D$.  Let
$\star$ be a semistar operation on $D$.  We can define a semistar
operation $\dot{\star}^{\mbox{\tiny \it \tiny T}}:
\overline{\boldsymbol{F}}(T)
\rightarrow \overline{\boldsymbol{F}}(T)$ on $T$ , by setting:

\centerline{ $E^{\dot{\star}^{\mbox{\tiny \it \tiny T}}} := E^\star\,, \; \mbox{
for each } \; E \in \overline{\boldsymbol{F}}(T) (\subseteq
\overline{\boldsymbol{F}}(D))\,; $}

\noindent \cite[Proposition 2.8]{FL01}.    When $T = D^\star$, we
denote simply by $\dot{\star}$ the (semi)star operation
$\dot{\star}^{\mbox{\tiny \it \tiny D$^\star$}}$ on $D^\star$.  \\
Note that ${(\dot{\star}_{_{\!f}})^{\mbox{\tiny \it \tiny T}}} =
({\dot{\star}^{\mbox{\tiny \it \tiny T}}})_{_{\!f}}$   \cite[Lemma
3.1]{FJS03}.  In particular,
if $\star = \star_{_{\!f}}$ then ${\dot{\star}^{\mbox{\tiny \it \tiny
T}}}$ is a semistar operation of finite type on $T$.

\bf (6) \rm On the other hand, if $\ast$ is a semistar operation on an
overring $T$ of an integral domain $D$, we can construct a semistar operation \
\d{$\ast$}$_{\mbox{\tiny \it \tiny D}}: \overline{\boldsymbol{F}}(D)
\rightarrow \overline{\boldsymbol{F}}(D)$ on $D$, by setting:

\centerline{$ E^{{\mbox{\d{$\ast$}}}_{\mbox{\tiny \it \tiny D}}} := (ET)^\ast\,,
\; \mbox{ for each } \; E \in \overline{\boldsymbol{F}}(D) \,, $}

\noindent\cite[Proposition 2.9]{FL01}. \\
  For more details on the semistar operations considered in (5) and (6),
cf.    \cite{Okabe/Matsuda:1994}  and \cite{FL01}.

\begin{rem} \label{new:rk} \rm  Let  $\star$  be a semistar operation on an integral 
domain $D$. For 
each nonzero ideal $I$ of $D$, define $I_{r(\star)}:= I^\star \cap D$. 
Then it is easy to see that the map\ $I \mapsto I_{r(\star)}$\  
defines a weak ideal system\ (= $x$--system in the sense of K. E. 
Aubert)\  on $D$ (as a commutative cancellative monoid, disregarding 
 the additive structure), cf. \cite[Chapter 2]{HK98}, therefore the 
 theory developed in \cite[Part A]{HK98} applies. In particular, 
 $r(\star_{_{\!f}}) = r(\star)_{s}$ \cite[page 25]{HK98},\ $\mathcal{M}(\star) 
 = r(\star)$--max$(D)$ \cite[page 57]{HK98}, and 
 $\widetilde{\star} =  r(\star)_{s}[d]$ \cite[Definition 3.1 and Proposition 3.2]{HK:2000}.
 
 Furthermore, using the more general setting of module systems 
 on monoids, the spectral semistar operations (Example \ref{ex:1.1} (3)) and the semistar operations\ ${\dot{\star}^{\mbox{\tiny \it \tiny 
 T}}}$\ and\ ${{\mbox{\d{$\ast$}}}_{\mbox{\tiny \it \tiny D}}}$,  defined 
 in Example \ref{ex:1,4-6} (5) and (6), have a natural correspondent 
 interpretation in terms of module systems, which is described in 
 \cite{Koch:JA}, and so the theory developed in this paper  also applies.
	\end{rem}

\begin{prop} \label{pr:1.2}
 Let $D$ be an integral domain and $T$ an overring of $D$.
\begin{enumerate}
\item[(1)] Let $\ast$ be a semistar operation on\  $T$.
 Denote simply by $\boldsymbol{\star}$ the semistar
operation \d{${\ast}$}$_{{\mbox{\tiny \it \tiny D}}}$, defined
on $D$, then the semistar operations \,
 $\dot{\boldsymbol{\star}}^{\mbox{\tiny \it T}}\,$ and $\, \ast \, $
 (both defined on $T\,$) coincide.
 \item[(2)] Let  $\star$   be a semistar operation on\ $D$.  Denote
simply by $\boldsymbol{\ast}$ the semistar operation
$\dot{\star}^{\mbox{\tiny \it T}}$, defined on $T$, then $\star \leq
\mbox{\d{$\boldsymbol{\ast}$}$_{{\mbox{\tiny \it \tiny D}}}$}$
(note that both
semistar operations are defined on $D$).
Furthermore, if $T = D^\star$ then  $\star =
\mbox{\d{$\boldsymbol{\ast}$}$_{{\mbox{\tiny \it \tiny D}}}$}$.
 \end{enumerate}
\end{prop}
   \begin{proof}  (1) and the first statement in (2) are
already in  \cite[Corollary 2.10]{FL01},
\cite[Lemma 45]{Okabe/Matsuda:1994}.  For the last statement
note that, for each $E \in  \overline{\boldsymbol{F}}(D)$,
$E^{\mbox{\d{$\boldsymbol{\ast}$}}_{{\mbox{\tiny \it \tiny D}}}} =
(ET)^\ast =(ET)^{\dot{\star}^{\mbox{\tiny \it T}}}= (ET)^\star =
(ED^\star)^\star = (ED)^\star = E^\star$.
\end{proof}
\end{exs}

If $R$ is a ring and $X$ an indeterminate over $R$, then the ring
$R(X):=\{f/g \, \vert \; f,g \in R[X] \mbox{ and }
\boldsymbol{c}(g)=R \}$ (where $\boldsymbol{c}(g)$ is the content
of the polynomial $g$) is called the \emph{ Nagata ring} of $R$
\cite[Proposition 33.1]{Gil92}.  A more general construction of a Nagata
ring associated to a semistar operation defined on an integral domain $D$
was considered in \cite{FL03} (cf.  also \cite{Ka89}  and \cite[Chapter 
20, Ex. 4]{HK98},   for the star case).

\begin{prop}  \label{prop:nagata}
   \rm \  \cite[Proposition 3.1,
Proposition 3.4, Corollary 3.5, Theorem 3.8]{FL03}.
\it  Let $\star$ be a semistar
    operation on an integral domain $D$ and set $N(\star) := N_D(\star):=\{h
    \in D[X] \, \vert \; h \neq 0 \mbox{ and }
    (\boldsymbol{c}(h))^\star=D^\star \}$.    Let $\tilde{\star}$
    (respectively, $\dot{\tilde{\star}}$) be the spectral semistar operation
    defined on $D$ (respectively, $D^{\tilde{\star}}$), introduced in Example
    \ref{ex:1,4-6}((4) and (5)).     Then:
\begin{enumerate}
\item[(1)] $N(\star)$ is a saturated multiplicative subset of
$D[X]$ and $N(\star)=N(\star_{\!_f})=D[X] \smallsetminus \bigcup
\{Q[X] \, \vert \; Q \in \mathcal{M}(\star_{\!_f}) \}$.
\end{enumerate}
\noindent  Set $\Na(D,\star):=   D[X]_{N(\star)}$ and call this
integral domain  \rm{the Nagata ring of $D$ with respect to the
semistar operation $\star$}.
\begin{enumerate}
    \item[(2)] $ \Max(\Na(D,\star))    = \{Q[X]_{N(\star)}
\, \vert \; Q \in \mathcal{M}(\star_{\!_f}) \}$ \it  and
$\mathcal{M}(\star_{\!_f})$ coincides with the canonical image in\
$\Spec(D)$ of\ $\Max \left(\Na(D,\star) \right)$\,.
\item[(3)] $
\Na(D,\star)    = \bigcap \{D_Q(X) \, \vert \; Q \in
\mathcal{M}(\star_{\!_f})\}$\,.
 \item[(4)] $E^{\tilde{\star}} =
E\mbox{\rm Na}(D, \star) \cap K\,, $ \it  for each $E \in
\overline{\boldsymbol{F}}(D) \,.  $
\item[(5)] $\mathcal{M}(\star_{f}) = \mathcal{M}(\tilde{\star})\,.$
\item[(6)] $\Na(D, \star) = \Na(D, {\tilde{\star}}) = \Na(D^{\tilde{\star}},\dot{\tilde{\star}})
\supseteq D^{\tilde{\star}}(X)\,.$
\hfill $\Box$
\end{enumerate}
\end{prop}

It is clear that
$\Na(D,\star)=\Na(D,\star_{\!_f})$ and, when $\star=d$ (the identity
(semi)star operation) on $D$, then $\Na(D, d)=D(X)$.

\medskip

 Given a semistar operation $\star$ on an integral domain $D$, we say
that    $\,\star\,$ is an \emph{e.a.b. (endlich arithmetisch brauchbar)
semistar operation} of $D$ if, for all $ E,F,G \in \boldsymbol{f}(D)$, \
$(EF)^\star \subseteq (EG)^\star$\ implies that\ $F^\star \subseteq
G^\star$ \cite[Definition 2.3 and Lemma 2.7]{FL01}.

   We recall next the definition of two relevant semistar operations,
 associated to a given semistar operation.

 \setcounter{thm}{0}

 \begin{exs} \label{ex:1,7-8} \bf (7) \rm  Given a
 semistar operation $\star$ on an integral domain $D$, we call \it the
 semistar integral closure $[\star]$ of $\star$, \rm the semistar
 operation on $D$ defined by setting:

  \centerline{$ F^{[\star]} := \cup\{((H^\star:H)F)^{\star_f}
  \,\;|\;\, H \in \boldsymbol{f}(D) \} \,,\;\;
 \mbox{for each } \, 
 F \in \boldsymbol{f}(D) \,,
$}

\noindent and

\centerline{$ E^{[\star]} := \cup \{F^{[\star]} \,\;|\;\, F \in
\boldsymbol{f}(D),\; F \subseteq E\}\,,\; \; \mbox{for each } \, E \in
\boldsymbol{\overline{F}}(D) \,.
$}

\noindent It is not difficult to see that the operation $[\star]$ is a semistar
operation of finite type on $D$, that $\, \star_{f} \leq [\star]\,$,
hence $\, D^{\star} \subseteq D^{[\star]}\,$, and that $D^{[\star]}$ is
integrally closed \cite[Definition 4.2, Proposition 4.3 and Proposition
4.5 (3)]{FL01}.  Therefore, it is obvious that if $\, D^{\star} =
D^{[\star]}\, $ then $D^{\star}$ is integrally closed.

\bf (8) \rm Given an arbitrary semistar operation $\star$ on an integral domain
$D$, it is possible to associate to $\star$, an e.a.b. semistar
operation of finite type $\, \star_a \, $ on $D$, called \it the e.a.b.
semistar operation associated to $\star$, \rm defined as follows:

\centerline{$ F^{\star_a} :=
 \cup\{((FH)^\star:H)\,\; |\; \, H \in \boldsymbol{f}(D) \}, \; \; \textrm {for each } \, 
 F \in \boldsymbol{f}(D) \, ,
$}

\noindent and

\centerline{$ E^{\star_a} := \cup\{F^{\star_a} \,\; |\; \, F \subseteq E\,,\; F
\in \boldsymbol{f}(D) \},\; \; \textrm {for each } \, E \in
\boldsymbol{\overline{F}}(D), $}

\noindent \cite [Definition 4.4]{FL01}.  Note that $\, [\star] \leq
{\star_a}\,$, that $D^{[\star]} = D^{\star_a}$ and if $\star$ is an
e.a.b. semistar operation of finite type then $\,\star = \star_{a}\,$
\cite [Proposition 4.5]{FL01}.\\
 More information about the semistar operations $[\star]$ and
${\star_a} $ can be found in  \cite{Jaffard}, \cite{OM1},
\cite{Okabe/Matsuda:1994}, \cite{Koch:1997}, \cite{HK98} ,
\cite{Koch:2000}, and \cite{FL01b}.
\end{exs}

Let $\star$ be a semistar operation on $D$ and
let $V$ be a valuation overring of $D$.  We say that $V$ is a \it $\star$--valuation overring of
$D$ \rm if, for each $F \in \boldsymbol{f}(D)\,$, $ F^\star \subseteq FV\,$ (or equivalently,
${\star_{_{\!f}}} \leq \star_{\{V\}}$ (Example~\ref{ex:1.1}(1)). \\
 Note that   a valuation overring $V$ of $D$ is a
$\star$--valuation overring of $D$ if and only if $V^{\star_{_{\!f}}} =V$.
  (The ``only if'' part is obvious; for the ``if'' part recall
 that, for each $F \in \boldsymbol{f}(D)$, there exists a nonzero
element $x \in K$ such that $FV = xV$, thus $F^\star \subseteq
(FV)^{\star_{_{\!f}}} =(xV)^{\star_{_{\!f}}} = xV^{\star_{_{\!f}}} =xV=
FV$).
\\ More details on semistar valuation overrings can be found in
\cite{FL01b}, \cite{FL03} (cf.  also \cite{Jaffard},  \cite{Koch:1997} and
\cite{Koch:2000}).

\medskip

We recall next the construction of the Kronecker function ring with 
respect to a semistar operation (the star case is studied in detail in 
\cite[Section 32]{Gil92} and  \cite[Chapter 20, Ex. 6]{HK98}).

\setcounter{thm}{4}
  \begin{prop}\label{prop:Kr} \rm
  \cite[Proposition 3.3, Theorem 3.11, Theorem 5.1, Corollary 5.2, Corollary
  5.3]{FL01}, \cite[Theorem 3.5]{FL01b}.
  \it Let $\star$ be any semistar
  operation defined on an integral domain $D$ with quotient field $K$ and
  let $\star_{a}$ be the e.a.b. semistar operation associated to $\star$
  (Example~\ref{ex:1,7-8}(8)).  Consider the e.a.b. (semi)star operation
  $\dot{\star}_{a} :={{\dot{\star}}_{a}}^{{\mbox{\tiny \it \tiny
  D}}^{\star_{a}}}$ (defined in Example~\ref{ex:1,4-6}(5)) on the integrally
  closed integral domain $D^{\star_{a}} = D^{[\star]}$ (Example~\ref{ex:1,7-8}((7) and (8))).  Set
\[
\begin {array} {rl}
\Kr(D,\star) := \{ f/g \; \,|\, & f,g \in D[X] \setminus \{0\} \;\;
\mbox{ \it and there exists } \; h \in D[X] \setminus \{0\} \; \\ &
\mbox{ \it such that } \; (\boldsymbol{c}(f)\boldsymbol{c}(h))^\star
\subseteq (\boldsymbol{c}(g)\boldsymbol{c}(h))^\star \,\} \, \cup\,
\{0\}\,.
\end{array}
  \]
Then we have:
\begin{enumerate}
\item[(1)]\ $\Kr(D,\star)$ is a B\'ezout domain with quotient field
$K(X)\,,$ called \rm the Kronecker function ring of $D$ with respect to
the semistar operation $\star\,.$ \it
\item[(2)]\  ${\Na}(D,\star) \subseteq {\Kr}(D,\star)\,.$
\item[(3)]\  ${\Kr}(D,\star) = {\Kr}(D,\star_a) = {\Kr}(D^{\star_a},\dot{\star}_{a})\,$.
\item[(4)]\
$ E^{\star_{a}} = E{\Kr}(D,\star) \cap K \,,$ \; for each $ E \in
\boldsymbol{\overline{F}}(D)$\,.
\item[(5)]\ $\Kr(D, \star) = \bigcap \{ V(X) \mid V \mbox{ is a $\
\star$--valuation overring of } D \}\,.$
\item[(6)]\ If $\,F := (a_{0},a_{1},\ldots, a_{n}) \in
 \boldsymbol{f}(D)$\, and $\,f(X) :=a_{0}+ a_{1}X +\ldots +a_{n}X^n \in
 K[X]\,,$ then:

 \centerline{ $ F{\Kr}(D,\star) = f(X){\Kr}(D,\star) =
 {\boldsymbol{c}}(f){\Kr}(D,\star)\,.  $}

 \end{enumerate} \vskip -0.5cm \hfill $\square$

  \end{prop} \rm

If $D$ is an integral domain and $\star$ is a semistar operation
on $D$, we say that a  nonzero    ideal $I$ is
\emph{$\star$--invertible}, if $(II^{-1})^\star= D^\star$.  We
define $D$ to be a \emph{P$\star$MD} if each  nonzero    finitely
generated  ideal of $D$ is
$\star_{_{\!f}}$--invertible  (cf.  \cite{FJS03} and also
\cite{Griffin:1967}, \cite{Mott/Zafrullah:1981},
\cite{Houston/Malik/Mott:1984}, \cite{Ka89}, \cite{GJS99}).  In 
particular, note that, if $\star$ is a star operation, then $D$ is a 
P$\star$MD if and only if D is $\star$--Pr\"{u}fer in the sense of 
\cite[Chapter 17]{HK98}.  

By
using $\Na(D,\star)$ and $\Kr(D,\star)$, we have the following
characterization of a P$\star$MD.

\begin{prop} \label{pr:1.5} \cite[Theorem 3.1, Remark 3.1]{FJS03}
Let $D$ be an integral domain and $\star$ a semistar operation on
$D$. The following are equivalent:
\begin{enumerate}
\item[(i)] $D$ is a P$\star$MD;
\item[(ii)] $D_Q$ is a valuation domain, for each $Q \in
\mathcal{M}(\star_{_{\!f}})$;
\item[(iii)] $\Na(D,\star)$ is a Pr\"ufer domain;
\item[(iv)] $\Na(D,\star)=\Kr(D,\star)$;
\item[(v)] $\tilde{\star}$ is an e.a.b. semistar operation;
\item[(vi)] $\star_{_{\!f}}$ is stable and e.a.b.
\end{enumerate}
In particular, $D$ is a P$\star$MD if and only if it is a
P$\tilde{\star}$MD.    Moreover, in a P$\star$MD, $\tilde{\star} =
\star_{_{\!f}}$. \hfill $\Box$
\end{prop}

Let $D$ be an integral domain and $T$
an overring of $D$.  Let $\star$ be a semistar operation on $D$ and
$\star^\prime$ a semistar operation on $T$.  Then,  we say that    $T$ is
\emph{$(\star, \star^\prime)$--linked to} $D$, if

\centerline{$F^\star=D^\star \; \Rightarrow\; (FT)^{\star^\prime}
=T^{\star^\prime}\,,   $}

\noindent for each nonzero finitely generated ideal $F$ of $D$.
Finally, recall that we say that    $T$ is \emph{$(\star,
\star^\prime)$--flat over} $D$ if it is $(\star, \star^\prime)$--linked to
$D$ and, in addition, $ D_{Q\cap D}=T_Q $, for each
quasi--$\star^{\prime}_{_{\!f}}$  --maximal    ideal $Q$ of $T$.  More details on these
notions can be found in \cite{EF} (cf.  also  \cite{KP} and  \cite{HK03}).


\section{Semistar almost Dedekind Domains}


Let $D$ be an integral domain and $\star$ a  semistar operation on
$D$. We say that $D$ is a \emph{semistar almost Dedekind domain}
 (for short,    a \emph{$\star$--ADD}) if $D_M$ is a rank-one discrete
valuation  domain  (for short, DVR),   for each quasi--$\star_{_{\!f}}$--maximal ideal
$M$ of $D$.

Note that,  by definition, $\star$--ADD $=$
$\star_{_{\!f}}$--ADD and that,    if $\star=d$\  (= the identity
(semi)star operation),   we obtain the classical notion of ``almost
Dedekind domain''  (for short, ADD)   as in \cite[Section
36]{Gil92}.   Note that, If $\star=v$,  the $v$--ADDs coincide    with
the $t$--almost Dedekind domains studied by Kang  \cite[Section
4]{Ka89}; more generally, if $\star$ is a star operation, then $D$ is 
a  $\star$--ADD if and only if $D$ is a $\star$--almost Dedekind 
domain in the sense of \cite[Chapter 23]{HK98}.   Note also that, a field has only the identity (semi)star
operation and thus a field is, by convention,  a trivial example of a ($d$--)ADD
(since, in this case, $\mathcal{M}(d) = \emptyset$).  \\
 An analogous notion of generalized
almost Dedekind domain was considered in the language of ideal
systems on commutative monoids in \cite[Chapter 23]{HK98}.

\begin{rem} \label{rk:2.1}
Let $\star_1, \star_2$ be two semistar operations on
$D$ such that  $({\star_1})_{_{\!f}} \le ({\star_2})_{_{\!f}}$.
If $D$ is a $\star_1$--ADD, then $D$ is a $\star_2$--ADD. In particular:

  --\    $D$  is a    ADD \; $\Rightarrow$ \; $D$ is a $\star$--ADD,
 for each semistar operation $\star$ on $D$;

  -- \   if
$\star$ is a  (semi)star operation on $D$  (so, $\star \le v$),    then:

\hskip 1cm $D$  is a    $\star$--ADD \; $\Rightarrow$\; $D$ is a $v$--ADD
(and, hence, $D$ is integrally closed).

\noindent Note that,  in general,    for a semistar operation $\star$, a
$\star$--ADD may be not integrally closed.  For instance, let $K$ be a
field and $T:=K[\![X]\!]  = K + M$, where $M :=XT$ is the maximal ideal of
 the discrete valuation domain   $T$.  Set $D:=R+M$, where $R$ is
a non integrally closed integral domain    with quotient field $K$
(hence, $D$ is not integrally closed   \cite[Proposition
2.2(10)]{F80}).   Take $\star:=\star_{\{T\}}$  on $D$.     Then, we
have  $\star=\star_{_{\!f}}$, $\dot{\star}^T   =d_T$  is the identity (semi)star
operation on $T$ and    $\mathcal{M}(\star_f)=\{M\}$
 (by \cite[Lemma 2.3(3)]{FL03})    and $D_M=T$  \cite[Proposition
1.9]{F80}.   So $D$ is a $\star$--ADD which is not integrally closed
 (hence, in particular, $D$ is not an ADD).
\end{rem}

\begin{prop} \label{prop:2.2} Let $D$ be an integral domain,  which is
not a field,     and $\star$ a semistar operation on $D$.  Then:
\begin{enumerate}
\item[(1)] $D$ is a $\star$--ADD
if and only if $D_P$ is a DVR, for  each
quasi--$\star_{_{\!f}}$--prime ideal $P$ of $D$.
\item[(2)] If  $D$ is a
$\star$--ADD, then $D$ is a P$\star$MD and each
quasi--$\star_{_{\!f}}$--prime  of $D$    is  a
quasi--$\star_{_{\!f}}$--maximal  of $D$.
\item[(3)] Let $T$
be an overring of $D$ and $\star'$ a semistar operation on $T$.  Assume
that $D \subseteq T$ is a $(\star, \star')$--linked extension.  If $D$
is a $\star$--ADD, then $T$ is a  $\star'$--ADD.   \item[(4)] If $D$
is a $\star$--ADD, then $D^\star$ is a $\dot{\star}$--ADD.
\end{enumerate}
\end{prop}
\begin{proof} (1) It follows  easily    from the fact that each
quasi--$\star_{_{\!f}}$--prime is contained in a
quasi--$\star_{_{\!f}}$--maximal  \cite[Lemma 2.3(1)]{FL03}.   \\
  (2) is a straightforward consequence of (1) and of Proposition
\ref{pr:1.5} ((i)$\Leftrightarrow$(ii)).
\\
(3) Let  $N   \in {\mathcal M}(\star'_{_{\!f}})$, then  $(N\cap
D)^{\star_{_{\!f}}}    \neq D^\star$  \cite[Proposition 3.2]{EF}.
Let  $M   \supseteq N\cap D$ be a quasi--$\star_{_{\!f}}$--maximal
ideal of $D$.  We have  $D_M \subseteq D_{N\cap D}\subseteq T_N$.
   So  $T_N=D_{N \cap D}=D_M$, because $D_M$    is a DVR
 (by assumption $D$ is a $\star$--ADD).     From this proof we deduce
 also    that $N\cap D  \ (= M)$ is a
quasi--$\star_{_{\!f}}$--maximal ideal  of $D$,    for each
quasi--$\star'_{_{\!f}}$--maximal ideal $N$ of $T$.  \\ (4) It follows
from \cite[Lemma 3.1(e)]{EF} and (3).
\end{proof}

\begin{rem} \label{rem:2.1} \bf  (1)  \rm We will show that, for a converse of
Proposition \ref{prop:2.2}(2), we  will    need additional
conditions (cf.  Theorem \ref{thm:2.7}((1)$\Leftrightarrow$(3), (4)).  \\
\bf (2) \rm  The converse of Proposition \ref{prop:2.2}(4) is not true in
general.  Indeed, let $K$ be a field and $k\subset K$ a proper subfield
of $K$.  Let $T :=   K[[X]]$ and $D :=   k+M$, where $M  :=
XT$   is the maximal ideal of $T$.   Take $\star:=\star_{\{T\}}$ on
$D$.     Note that  $\star=\star_{_{\!f}}$ and that
$\dot{\star}^T=d_T$ is the identity (semi)star operation on $T$.  We
have that $T =D^\star$   is a  $\dot{\star}^T$--ADD $=$   ADD
(since $T$ is a DVR), but $D$ is not a $\star$--ADD, since $M$
is a  quasi--$\star_{_{\!f}}$--maximal    ideal  of $D$  and  (by
\cite[Proposition 1.9]{F80})    $D_M=D$ is not a valuation domain.
\end{rem}

\begin{prop} \label{prop:2.4} Let $D$ be an integral domain and $\star$ a \sl
(semi)star \it operation on $D$.  Then the following are equivalent:
\begin{enumerate}
\item[(1)] $D$ is a $\star$--ADD.
\item[(2)] $D$ is a $t$-ADD and $\star_f=t$.
\end{enumerate}
\end{prop}
\begin{proof} (1) $\Rightarrow$ (2) By Remark \ref{rk:2.1}, if $D$ is a
$\star$--ADD, then $D$ is a $v$--ADD or, equivalently, a $t$--ADD.
Moreover, by Proposition \ref {prop:2.2}(2) and \cite[Proposition
3.4]{FJS03}, $\star_f=t$.  The converse is clear.
\end{proof}

%

\begin{thm} \label{thm:2.3} Let $D$ be an integral domain,  which is
not a field,    and $\star$ a semistar operation on $D$.  The following are
equivalent:
\begin{enumerate}
\item[(1)] $D$ is a $\star$--ADD.
\item[(2)] $\Na(D,\star)$ is an ADD (i.e. $\Na(D,\star)$ is a
1-dimensional Pr\"ufer domain and contains no idempotent maximal
ideals).
\item[(3)] $\Na(D,\star) = \Kr(D, \star)$ is an ADD and a B\'ezout
domain.
\end{enumerate}
\end{thm}
\begin{proof} (1) $\Leftrightarrow$ (2).  By   Proposition
\ref{prop:nagata}(2),
    the    maximal ideals of $\Na(D,\star)$ are of the form $M\!
\Na(D,\star)$, where $M\in{\mathcal M}(\star_{_{\!f}})$.  Also, for each
$M\in{\mathcal M}(\star_{_{\!f}})$, we have $\Na(D,\star)_{M\!\Na(D,
\star)}=D_M(X)$.  Moreover, it is well-known that, for $M\in{\mathcal
M}(\star_{_{\!f}})$, $D_M$ is a DVR if and only if $D_M(X)$ is a DVR
 \cite [Theorem 19.16 (c), Proposition 33.1 and Theorem 33.4
((1)$\Leftrightarrow$(3))]{Gil92}.     From  these facts    we
conclude  easily.   \\
(1)$\Rightarrow$(3). If $D$ is a $\star$--ADD,  in particular $D$ is
a P$\star$MD  (Proposition \ref{prop:2.2}(2)),    then
$\Na(D,\star)=\Kr(D,\star)$, by   Proposition \ref{pr:1.5} ((i )$\Leftrightarrow$(iv)).
Therefore, we deduce that  $\Na(D,\star)$
is a B\'ezout domain (Proposition \ref{prop:Kr}(1)) and an ADD by
(1)$\Rightarrow$(2).   \\
(3) $\Rightarrow$ (2) is trivial.
\end{proof}

\begin{cor} \label{cor:2.6} Let $D$ be an integral domain and $\star$ a semistar
operation on $D$.     The following are equivalent:
\begin{enumerate}
\item[(1)] $D$ is a $\star$--ADD.
\item[(2)] $D$ is a $\tilde{\star}$--ADD.
\item[(3)] $D^{\tilde{\star}}$ is a $\dot{\tilde{\star}}$--ADD.
 \item[(4)]  $D^{\tilde{\star}}$ is a $t$--ADD and $\dot{\tilde{\star}}=t_{
D^{\tilde{\star}}}$.
\end{enumerate}
\end{cor}
\begin{proof} Note that $\Na(D, \star)=\Na(D, \tilde{\star})=
\Na(D^{\tilde{\star}}, \dot{\tilde{\star}})$   (Proposition
\ref{prop:nagata}(6)),
then apply Theorem \ref{thm:2.3}((1)$\Leftrightarrow$(2))
 to obtain the equivalences $(1)\Leftrightarrow(2)\Leftrightarrow(3)$.
The equivalence (3)$\Leftrightarrow$(4) follows from Proposition
\ref{prop:2.4}.
\end{proof}

 Next goal is a characterization of $\star$--ADD's in terms of
valuation overrings,  in the style of \cite[Theorem 36.2]{Gil92}.  For
this purpose, we prove preliminarly the following:

\begin{lemma} \label{lemma:2.2}
Let $D$ be an integral domain and $\star$ a semistar operation on
$D$. Let $V$ be a valuation overring of $D$. Then the following
are equivalent:
\begin{enumerate}
\item[(1)] $V$ is a $\tilde{\star}$--valuation overring of $D$.
\item[(2)] $V$ is $(\tilde{\star}, d_V)$--linked to $D$.
\end{enumerate}
\end{lemma}
\begin{proof}
(1) $\Rightarrow$ (2): Since $V$ is a $\tilde{\star}$--valuation
overring, then $\tilde{\star} \leq \star_{\{V\}}$.  Thus,  the present
implication    follows from the fact that
$\dot{\star}^{V}_{\{V\}}   = d_V$  (so $\dot{\tilde{\star}}^{V} =
d_V$)    and  from \cite[Lemma 3.1(e)]{EF}.     \\
(2) $\Rightarrow$ (1): Let  $N$    be the maximal ideal of $V$
(which is $(\tilde{\star}, d_V)$--linked to $D$).     Then  $(N \cap
D)^{\tilde{\star}} \neq D^{\tilde{\star}}$   by \cite[Proposition
3.2  ((i)$\Rightarrow$(v))  ]{EF}.  Thus, there exists  $M    \in
\mathcal{M}(\star_{_{\!f}})  = \mathcal{M}(\tilde{\star})$
 (Proposition \ref{prop:nagata}(5))
such that  $N \cap D \subseteq M$.
   Hence  $D_M \subseteq D_{N \cap D}    \subseteq V$.  So, if $F
\in \boldsymbol{f}(D)$, then $F^{\tilde{\star}} \subseteq  FD_M
\subseteq FV$.  Therefore, $V$ is a $\tilde{\star}$--valuation overring
 of $D$.
\end{proof}

\begin{thm} \label{thm:2.2}
Let $D$ be an  integral    domain,  which is not a field,     and
$\star$ a semistar operation on $D$.  The following are equivalent:
\begin{enumerate}
\item[(1)] $D$ is $\star$--ADD.
\item[(2)] $D^{\tilde{\star}}$ is integrally closed and each $\tilde{\star}$--valuation overring of $D$ is a
 DVR.
\item[(3)] $D^{\tilde{\star}}$ is integrally closed and each
valuation overring $V$  of $D$, which is    $(\tilde{\star},d_V)$--linked to
$D$, is a DVR. \item[(4)] $D^{\tilde{\star}}$ is integrally
closed and each valuation overring $V$  of $D$, which is    $(\star,
d_V)$--linked to $D$, is a DVR.
\end{enumerate}
\end{thm}

\begin{proof} (1) $\Rightarrow$ (2). Since $D^{\tilde{\star}}=
\bigcap \{D_M \mid \, M \in \mathcal{M}(\star_{_{\!f}})\}$ and $D_{M}$
is a DVR, for each $ M \in \mathcal{M}(\star_{_{\!f}})$, then
$D^{\tilde{\star}}$ is integrally closed.     Now, let $V$ be a
$\tilde{\star}$--valuation overring of $D$, then $V\supseteq D_M$ for
some $M\in{\mathcal M}(\star_{_{\!f}})$ \cite[Theorem 3.9]{FL03}.  Since
$D_M$ is a DVR, then $V=D_M$ (is a DVR).\\
(2) $\Leftrightarrow$ (3). Follows immediately from Lemma
\ref{lemma:2.2}. \\
(3) $\Rightarrow$ (4).  It is  an  immediate consequence of the fact
that $\tilde{\star} \leq \star$ (cf.  \cite[Lemma 3.1(h)]{EF}).\\
(4) $\Rightarrow$ (1). Let $M\in{\mathcal M}(\star_{_{\!f}})$ and $V$ be
valuation overring of $D_M$. Then $V=V_{D \smallsetminus M}$ is
$(\star, d_V)$--linked to $D$  (cf.  \cite[Example 3.4(1)]{EF}).  Hence,
 by assumption,    $V$ is a DVR. Furthermore, $D_M$ is integrally
closed, since  ${D^{\tilde{\star}}} \subseteq D_{M}$ and thus
$D_M = {D^{\tilde{\star}}}_{MD_M \cap D^{\tilde{\star}}}$.  So $D_M$ is
an ADD, by \cite[Theorem 36.2]{Gil92}, that is,  $D_M$    is a
 DVR. Therefore $D$ is a $\star$--ADD.
\end{proof}

\begin{cor} Let $D$ be an integral domain, which is not a field.  Then the
following are equivalent:
\begin{enumerate}
\item[(1)]$D$ is $t$--almost Dedekind domain.
\item[(2)] $D$ is integrally closed and each $w$--valuation overring of
$D$ is a DVR.
\item[(3)] $D$ is integrally closed and each $t$--linked valuation overring
of $D$ is a DVR.
\end{enumerate}
\end{cor}
\begin{proof} This is an immediate consequence of Theorem \ref {thm:2.2}
and of the  wellknown  fact that for a valuation domain $V$,
$d_V= w_{V}= t_V$  (cf. also \cite[Section 3]{EF} for the
$t$--linkedness).
\end{proof}

\begin{rem} If $D$ is a $\star$--ADD,  which is not a field,   then, by
Theorem \ref{thm:2.2}
    and by the fact that a $\star$--valuation overring is a
    $\tilde{\star}$--valuation overring, each $\star$--valuation overring of
    $D$ is a DVR. Note that the converse is not true, even if
    $D^{\tilde{\star}}$ is integrally closed.  Let $D$ and $T$ be as in
    Remark \ref{rem:2.1}(2).  Assume that $k$ is algebraically closed in
    $K$.  Since $\star = \star_{\{T\}}$, then $\star = \star_{_{\!f}}$,
    $\mathcal{M}(\star_{_{\!f}})=\{M\}$ and $D=D_M=D^{\tilde{\star}}$ is
    integrally closed, where $\tilde{\star}=d_D$.  Moreover,  each
    $\star$--valuation overring of $D$ is necessarily a valuation overring
    of $T$  (since $T = D^{\star_{_{\!f}}} =D^\star \subseteq V
    =V^{\star_{_{\!f}}} =V^\star$).  This implies that each $\star$--valuation overring of $D$ is a
     DVR  (since the only non trivial valuation overring of $T$
    is $T$, which is a    DVR).   Therefore, by
    Proposition \ref{prop:nagata}(6) and \ref{prop:Kr}(5),  
   $\Na(D,\star)=\Na(D^{\tilde{\star}},\dot{\tilde{\star}})=\Na(D,d_D)=
      D(Z)  \subsetneq \Kr(D,\star)=\Kr(T,d_T)=   T(Z) $  (where
     $Z$   is
    an indeterminate over $T$ and $D$).   On the other hand,    since
    $\mbox{t.deg}_k(K) \geq 1$, it is possible to find
      ($\tilde{\star}$--)   valuation overrings of $D$ (of rank $\geq 2$)
    contained in $T$  \cite[Theorem 20.7]{Gil92}.
\end{rem}

Let  $D$     be an integral domain and $\star$ a semistar operation on $D$.
For each quasi-$\star$-prime $P$ of $D$, we define the
\emph{$\star$-height} of $P$ (for short, $\star\mbox{-ht}(P)$) the
 supremum of the lengths of the chains  of quasi--$\star$--prime
ideals of $D$,  between prime ideal $(0)$ (included) and $P$.
Obviously, if $\star = d$ is the identity (semi)star operation on $D$,
then $d$-ht$(P) = $ ht$(P)$, for each prime ideal $P$ of $D$.  If the
set of  quasi--$\star$--primes  of $D$ is not empty,    the  \it
$\star$-dimension of $D$ \rm is defined as follows:

 \centerline{$\star \mbox{-dim}(D) := \Sup \{\star \mbox{-ht}(P) \mid
\;  P \, \mbox{ is a quasi--}\star\mbox{--prime of }D\}$\,.}

 \noindent If the set of  quasi--$\star$--primes   of $D$ is empty, then we
pose $\star \mbox{-dim}(D) := 0$.

\noindent Note that, if $\star_1 \leq \star_2$,  then $\star_2$-dim$(D) \leq
\star_1$-dim$(D)$.   In particular, $\star$-dim$(D) \leq d$-dim$(D) = $
dim$(D)$\ (= Krull dimension of $D$), for each semistar operation $\star$
on $D$.  Note that, recently, the notions of $t$-dimension and of
$w$-dimension have been received a considerable interest by several
authors (cf.  for instance,  \cite{H94},  \cite{W99} and \cite{W01}).

\begin{lemma} \label{lemma:2.5} Let $D$ be an integral domain and
$\star$ a semistar operation on $D$, then
$$
\begin{array}{rl}
    \tilde{\star}\mbox{-}\di(D)=& \hskip -5pt \Sup\{\htt(M) \mid \, M
    \in{\mathcal M}(\star_{_{\!f}}) = {\mathcal M}({\tilde{\star}})\}=\\
    =& \hskip -5pt  \Sup\{\htt(P) \mid \, P\, \mbox{ is a
    quasi--}{\tilde{\star}}\mbox{--prime of }D\}\,.
\end{array}$$
\end{lemma}
\begin{proof} Let $M \in{\mathcal M}(\star_{_{\!f}})$ and $P\subseteq M$
be a  nonzero  prime
ideal of $D$. Since
$\mathcal{M}(\star_{_{\!f}})=\mathcal{M}(\tilde{\star})$   (Proposition
\ref{prop:nagata}(5))
we have $P \subseteq P^{\tilde{\star}}\cap D \subseteq
PD_M \cap D  =    P $.  So $P$ is a quasi--$\tilde{\star}$--prime ideal of
$D$.  Hence $\htt(M) =  \tilde{\star}\mbox{-}\htt(M)$,     so we get the
Lemma.
\end{proof}

\begin{rem}  Note that, in general,

    \centerline{$\star_{_{\!f}} \mbox{-dim}(D) \leq \Sup \{\htt(P)
    \mid \; P \, \mbox{ is a quasi--}\star_{_{\!f}}\mbox{--prime of
    }D\}$\,.}

\noindent    Moreover, it can happen that $\star_{_{\!f}} \mbox{-dim}(D)
\lneq \Sup \{\htt(P) \mid \; P \, \mbox{ is a
quasi--}\star_{_{\!f}}\mbox{--prime}$ of $D\}$,  as shows the
following example.\\
Let $T$ be a DVR , with maximal ideal $N$, dominating a
two-dimensional local Noetherian domain $D$, with maximal ideal $M$
\cite{C54} (or \cite[Theorem]{CHL96}), and let $\star := \star_{\{T\}}$.  Then,
clearly, $\star = \star_{_{\!f}}$ and the only
quasi--$\star_{_{\!f}}$--prime ideal of $D$ is $M$, since if $P$ is a
nonzero prime ideal of $D$, then $P^\star =PT =N^k$, for some integer
$k\geq 1$.  Therefore, if we assume that $P$ is
quasi--$\star_{_{\!f}}$--ideal of $D$, then we would have $P = PT \cap D
=N^k \cap D \supseteq M^k$, which implies that $P = M$.  Therefore, in
this case, $ 1= \star_{_{\!f}} \mbox{-dim}(D) =\star_{_{\!f}}$-ht$(M)
\lneq \Sup \{\htt(P) \mid \; P \, \mbox{ is a
quasi--}\star_{_{\!f}}\mbox{--prime}$ of $D\} = $ ht$(M) = \di(D) =2$.
Note that, in the present example, $\tilde{\star} $ coincides with the
identity (semi)star operation on $D$.  \\
 It is already known that, when $\star =v$, it may happen that
$t$--$\di(D) < w$--$\di(D)$, \cite[Remark 2]{W01}.

\end{rem}

The following  lemma    generalizes \cite[Theorem 23.3,  the first
statement in    (a)]{Gil92}.

\begin{lemma} \label{lemma:2.6} Let $D$ be a P$\star$MD. Let $Q$ be a
 nonzero     $P$--primary ideal of $D$,  for some prime ideal $P$ of
$D$, and let    $x\in D\smallsetminus P$.  Then
$Q^{\tilde{\star}}={(Q(Q+xD))}^{\tilde{\star}}$.
\end{lemma}
\begin{proof} Let $M\in \mathcal{M}(\star_{_{\!f}})$.  If $Q
\not\subseteq    M$, then $QD_M=Q^2D_M=Q(Q+xD)D_M(=D_M)$.  If
$Q\subseteq M$, then $QD_M$ is $PD_M$-primary and $x\in
D_M\smallsetminus PD_M$; so $QD_M=QxD_M$, by \cite[Theorem
17.3(a)]{Gil92},  since $D_M$ is a valuation domain.  Thus
$QD_M=(Q^2+Qx)D_M$,
 hence 
$Q^{\tilde{\star}}={(Q(Q+xD))}^{\tilde{\star}}$.
\end{proof}

 Let $\star$ be a semistar operation on an integral domain $D$. We
say that $D$ has  \it the    ${\star}$--cancellation law \rm  (for
short, \it
${\star}$--CL\rm ) if  $A, B, C\in
\boldsymbol{F}(D)$ and
$(AB)^{{\star}}=(AC)^{{\star}}$ implies that
$B^{{\star}}=C^{{\star}}$.  The following  theorem provides several characterizations of the
semistar almost Dedekind domains and,   in particular,    it  generalizes \cite[Theorem
36.5]{Gil92} and \cite[ Theorem 4.5]{Ka89}.

\begin{thm} \label{thm:2.7} Let $D$ be an integral domain  which is
not a field    and let $\star$ be a semistar operation on $D$.  The
following are equivalent:
\begin{enumerate}
\item[(1)] $D$ is $\star$-ADD.
\item[(2)] $D$ has  the    $\tilde{\star}$--cancellation law.
\item[(3)] $D$ is a P$\star$MD,
$\star_{_{\!f}}\mbox{-}\di(D)=1$ and $(M^2)^{\star_{_{\!f}}}\neq
M^{\star_{_{\!f}}}$, for each $M\in \mathcal{M}(\star_{_{\!f}})
\   (=\mathcal{M}(\tilde{\star}))$.
\item[(4)] $D$ is a P$\star$MD and
$\cap_{n\geq 1}(I^n)^{\star_{_{\!f}}}=0$ for each proper
quasi--$\star_{_{\!f}}$--ideal $I$ of $D$.
\item[(5)] $D$ is a P$\star$MD and it has  the    ${\star}_{_{\!f}}$--cancellation law.
\end{enumerate}
\end{thm}
\begin{proof}
(1) $\Rightarrow$ (2). Let $A, B, C$ be three nonzero (fractional)
ideals of $D$ such that
$(AB)^{\tilde{\star}}=(AC)^{\tilde{\star}}$. Let $M\in
\mathcal{M}(\star_{_{\!f}})$.  Then,  we have    $ABD_M=
(AB)^{\tilde{\star}}D_{M}= (AC)^{\tilde{\star}}D_{M} =    ACD_M$  (we used
twice the fact that $\tilde{\star}$ is spectral, defined by
$\mathcal{M}(\star_{_{\!f}})$).
Moreover,    since $D_M$ is a DVR then,  in particular,    $AD_M$
is principal, thus $BD_M=CD_M$.  Hence
$B^{\tilde{\star}}=C^{\tilde{\star}}$.\\
(2) $\Rightarrow$ (3). If $D$ has $\tilde{\star}$--CL, then in
particular, $\tilde{\star}$ is an e.a.b. semistar operation on
$D$ \cite[Lemma 2.7]{FL01}, thus $D$ is a P$\star$MD
  (Proposition \ref{pr:1.5} ((v)$\Rightarrow$(i))\ ).
Let
$M\in\mathcal{M}(\star_{_{\!f}})$.  Clearly,  by $\tilde{\star}$--CL,
   $(M^2)^{\tilde{\star}}\neq M^{\tilde{\star}}$, and hence
$(M^2)^{\star_{_{\!f}}}\neq M^{\star_{_{\!f}}}$ (since
$\tilde{\star}=\star_{_{\!f}}$   by  Proposition \ref{pr:1.5}).
Next we show that ht$(M) =
1$, for each $M \in \mathcal{M}(\star_{_{\!f}})$.  Deny, let $P\subset
M$ be a nonzero prime ideal of $D$ and let $x\in M\smallsetminus P$.  By
Lemma \ref{lemma:2.6}, $P^{\tilde{\star}}={(P(P+xD))}^{\tilde{\star}}$.
Hence $D^{\tilde{\star}}={(P+xD)}^{\tilde{\star}}$, by
$\tilde{\star}$--CL. So $P+xD \not\subseteq M$, which is impossible.
Hence ht$(M)=1$, for each $M \in \mathcal{M}(\star_{_{\!f}})$.
Therefore, we conclude that $\star_{_{\!f}} \mbox{-dim}(D)=
{\tilde{\star}}\mbox{-dim}(D)=1$ (Lemma \ref{lemma:2.5}).\\
(3) $\Rightarrow$ (4).   Recall that each proper
quasi--$\star_{_{\!f}}$--ideal is contained in a
quasi--$\star_{_{\!f}}$--maximal ideal, then    it suffices to show that
$\cap_{n\ge 1}(M^n)^{\star_{_{\!f}}}=0$, for each $M\in
\mathcal{M}(\star_{_{\!f}})$.  Since,  by assumption
$(M^2)^{\star_{_{\!f}}}\neq M^{\star_{_{\!f}}}$,  then in particular
   $(M^2)^{\tilde{\star}}\neq M^{\tilde{\star}}$,  and so
$M^2D_M\neq MD_M$.  Henceforth $\{M^nD_M\}_{n\ge 1}$ is the set of
$MD_M$-primary ideals of $D_M$ \cite[Theorem 17.3(b)]{Gil92}.
 From the assumption we deduce that    dim$(D_M)=1$ (because
$\star_{_{\!f}}=\tilde{\star}$   by
 Proposition \ref{pr:1.5}),
then $\cap_{n\geq 1}M^nD_M=0$
\cite[Theorem 17.3 (c) and (d)]{Gil92}.  In particular, we have
$\cap_{n\ge 1}(M^n)^{\tilde{\star}}\subseteq \cap_{n\ge
1}\left((M^n)^{\tilde{\star}}D_M\right) = \cap_{n\ge
1}\left(M^nD_M\right) =0$, therefore $\cap_{n\ge
1}(M^n)^{\star_{_{\!f}}}=0$.  \\
(4) $\Rightarrow$ (1).  Let $M\in \mathcal {M}(\star_{_{\!f}})$.   It
is easy to see that     $(M^n)^{\tilde{\star}}= M^nD_M\cap D^{\tilde{\star}}$,
for each $n\ge 1$.  So, $(\cap_{n\ge 1}M^nD_M)\cap D^{\tilde{\star}}
\subseteq     \cap_{n\ge 1}(M^nD_M\cap D^{\tilde{\star}})=\cap_{n\ge
1}(M^n)^{\tilde{\star}} \subseteq \cap_{n\ge 1}(M^n)^{\star_{_{\!f}}}=0$
(the last equality holds by assumption).     Hence $\cap_{n\ge 1}M^nD_M=0$,
since $D_M$ is  an essential valuation    overring of
$D^{\tilde{\star}}$.  It follows that $D_M$ is a DVR \cite[p.
192 and Theorem 17.3(b)]{Gil92}. \\
  (2) $\Leftrightarrow$ (5) is a consequence of the fact that in a
P$\star$MD, $\widetilde{\star} = \star_{_{\!f}}$ and that the
$\tilde{\star}$--CL  implies P$\star$MD.
\end{proof}

\begin{rem}
  As a comment to Theorem \ref{thm:2.7} ((1)$\Leftrightarrow$(5)),
note that
$D$  may have  the $\star_{_{\!f}}$--CL without
being a $\star$--ADD.    It is sufficient to consider the example in Remark
 \ref{rem:2.1}(2).     In that case, $\star=\star_{_{\!f}}$ and
$\tilde{\star}=d_D$, since $\mathcal{M}(\star_{_{\!f}})=\{M\}$.
Clearly, $D$ has the $\star$--cancellation law (because $T$ is a
 DVR), but, as we have already remarked, $D$ is not a
$\star$--ADD, hence,  equivalently,    $D$ has not the
($\tilde{\star}$--)cancellation law.
\end{rem}

 Next result provides a generalization to the semistar case of
\cite[Theorem 36.4 and Proposition 36.6]{Gil92}.

\begin{prop} \label{prop:2.8} Let $D$ be an integral domain,  which
is not a field,    and $\star$ a semistar operation on $D$.  The following are
equivalent:

\begin{enumerate}
\item[(1)] $D$ is a $\star$--ADD.
\item[(2)] For each nonzero ideal $I$  of $D$, such that
$I^{\star_{_{\!f}}}\neq D^\star$ and $\sqrt{I}=:P$ is a prime ideal of
$D$,     then $I^{\tilde{\star}}=(P^n)^{\tilde{\star}}$, for some $n\ge 1$.
\item[(3)] $\tilde{\star}\mbox{-}\di(D)=1$ and, for each primary
quasi--$\tilde{\star}$--ideal $Q$  of $D$, then
$Q^{\tilde{\star}}=(M^n)^{\tilde{\star}}$, for some $M \in
\mathcal{M}(\star_{_{\!f}})$ and for some $n\ge 1$.
\end{enumerate}
\end{prop}
\begin{proof}
(1) $\Rightarrow$ (2) and (3).  Let  $I$ be a nonzero ideal  of $D$ with
$I^{\star_{_{\!f}}}\neq D^\star$ and $\sqrt{I}=P$ is prime.    Let $M$ be a
quasi--$\star_{_{\!f}}$--maximal ideal of $D$ such that $I\subseteq M$.
So $\sqrt{I}=P\subseteq M$, and hence $P=M$, since $D_M$ is a 
DVR. Thus $ID_M=M^nD_M$ for some $n\ge 1$.  On the other hand, if $N\in
\mathcal{M}(\star_{_{\!f}})$ and $N\neq M$, then $ID_N=D_N=M^nD_N$.
Hence $I^{\tilde{\star}}=(M^n)^{\tilde{\star}}$,  i.e.
$I^{\tilde{\star}}= (P^n)^{\tilde{\star}}$.     The fact that
$\tilde{\star}\mbox{-dim}(D)=1$ follows from
Theorem \ref{thm:2.7}((1)$\Rightarrow$(3))  (since, in the present situation,
$\star_{_{\!f}}=\tilde{\star}$).    \\
(2) $\Rightarrow$ (1). Let  $M\in \mathcal{M}(\star_{_{\!f}})$. Let $A$ be
an ideal of $D_M$ and assume that $\sqrt{A}=PD_M$, for some prime
ideal $P$ of $D$, $P \subseteq M$.
 Set $B:=A\cap D$. We have  $\sqrt{B}=P$ and hence
 $B^{\star_{_{\!f}}}\subseteq M^{\star_{_{\!f}}}\subset D^\star$. By
assumption, $B^{\tilde{\star}}=(P^n)^{\tilde{\star}}$, for some
$n\ge 1$,
 hence $A=(A\cap D)D_M=BD_M=B^{\tilde{\star}}D_M=(P^n)^{\tilde{\star}}D_M=P^nD_M$.
It follows from \cite[Proposition 36.6]{Gil92} that $D_M$ is an
ADD. Hence
$D_M$ is a DVR. \\
 (3) $\Rightarrow$ (1). We can assume $\star=\star_{_{\!f}}$, since  $\star$--ADD
 and $\star_{_{\!f}}$--ADD coincide. Let
   $M\in \mathcal{M}(\star_{_{\!f}})\   ( = \mathcal{M}(\tilde{\star})$
     (Proposition \ref{prop:nagata}(5)).
   Since
   $\tilde{\star}\mbox{-dim}(D)=1$, then   $\htt(M) =    \di(D_M)=1$
   (Lemma \ref{lemma:2.5}).  We can now proceed and conclude as in the
   proof of (2) $\Rightarrow$ (1).  (In this case, we  have
   $\sqrt{A}=MD_M$ and so $B$ is a $M$--primary
   quasi--$\tilde{\star}$--ideal  of $D$.   Therefore,  by
   assumption,    $B^{\tilde{\star}}=(M^n)^{\tilde{\star}}$, for some $n
   \geq 1$.)
\end{proof}

\begin{rem}
Note that, if $D$ is a $\star$--ADD,  which is not a field,  then
necessarily $D$   satisfies   the following conditions (obtained from the
statements (2) and (3) of Proposition \ref{prop:2.8}; recall that, in
this case, $\star_{_{\!f}}=\tilde{\star}$,  by Proposition \ref{prop:2.2}(2)
and
  Proposition \ref{pr:1.5}):
\begin{enumerate}
\item[(2$_{_{\!f}}$)]   For each nonzero ideal $I$   of $D$, such that
$I^{\star_{_{\!f}}}\neq D^\star$ and $\sqrt{I}=:P$ is a prime ideal of
$D$, then $I^{\star_{_{\!f}}}=(P^n)^{\star_{_{\!f}}}$, for some $n \ge
1$.
\item[(3$_{_{\!f}}$)] ${\star_{_{\!f}}}\mbox{-}\di(D)=1$ and, for
each primary quasi--${\star_{_{\!f}}}$--ideal $Q$ of $D$, then
$Q^{\star_{_{\!f}}}=(M^n)^{\star_{_{\!f}}}$, for some $M \in
\mathcal{M}(\star_{_{\!f}})$ and for some $n \ge 1$.
\end{enumerate}
  On the other hand, $D$ may satisfy either (2$_{_{\!f}}$) or
(3$_{_{\!f}}$) without being a  $\star$--ADD. 
It is sufficient to consider the example in Remark \ref{rem:2.1}(2).  In
that case, $\star =\star_{_{\!f}}$ and
$\mathcal{M}(\star_{_{\!f}})=\{M\}$.  Clearly, since $D$ is a local
one-dimensional domain (in fact,
$\tilde{\star}\mbox{-}\di(D)= \star_{_{\!f}}\mbox{-}\di(D)=\di(D) =1$), for each
nonzero ideal $I$  of $D$, with $I^{\star_{_{\!f}}}\neq D^\star$,
then $\sqrt{I}=M$ and $I^{\star_{_{\!f}}}=(M^n)^{\star_{_{\!f}}}$, for
some $n \ge 1$,  since $T$ is a DVR.   But, as we have already remarked, $D$ is not a
$\star$--ADD.
\end{rem}


\section{Semistar Noetherian domains}


Let $D$ be an integral domain and $\star$ a  semistar operation on
$D$. We say that $D$ is a \emph{$\star$--Noetherian domain} if  $D$
has the   ascending chain condition    on quasi--$\star$--ideals.

Note that, if $d\ (= d_D)$ is the identity (semi)star operation on $D$,
 the   $d$--Noetherian domains are just the usual Noetherian domains
and the notions of $v$--Noetherian  [respectively,
$w$--Noetherian]  domain and Mori  [respectively, strong Mori]
 domain coincide
\cite[Theorem 2.1]{Ba2000}  [respectively, \cite{WMc97}]. 

Recall that the concept of \sl star \rm Noetherian domain has already been
introduced, see for instance \cite{AA},  \cite{Z89} and \cite{GJS99}.    Using
ideal systems on commutative monoids, a similar general notion of
noetherianity  was considered
in \cite[Chapter 3]{HK98}.

\begin{lemma} \label{lemma:noeth} Let $D$ be an integral domain.
\begin{enumerate}
\item[(1)] Let $\star \le \star'$ be two semistar operations on $D$, then $D$
is $\star$--Noetherian  implies $D$ is  $\star'$--Noetherian. \\
 In particular:
\begin{enumerate}
\item[(1a)]  A Noetherian domain is a $\star$--Noetherian domain, for any
 semistar operation $\star$ on $D$.
\item[(1b)]  If $\star$ is a (semi)star operation and  if $D$ is a
$\star$--Noetherian domain, then $D$ is a Mori domain.
\end{enumerate}
\item[(2)] Let $T$ be an overring of $D$ and ${\ast}$ a semistar
operation on $T$.  If $T$ is   $\ast$  --Noetherian, then $D$ is
\rm{\d{${\ast}$}}${_{_{\!D}}}$--Noetherian.  In particular, if $\star$ is a semistar
operation on $D$, such that $D^\star$ is a $\dot{\star}$--Noetherian
domain, then $D$ is a $\star$--Noetherian domain.
\end{enumerate}
\end{lemma}
\begin{proof}
(1)  The first statement holds    because each quasi--$\star'$--ideal is a
quasi--$\star$--ideal.  (1a) and (1b) follow from (1)  since, for
each semistar operation $\star$, $d \leq \star$ and, if $\star$ is a (semi)star
operation, then $\star \leq v$.  \\
(2) If we have a chain of quasi--\d{${\ast}$}--ideals  $\{{I_n}\}_{n\geq 1}$
of $D$    that does not stop  then, by considering
$\{{(I_nT)^\ast}\}_{n\geq 1}$,    we get a chain of
quasi--${\ast}$--ideals  of $T$
that does not stop, since two distinct quasi--\d{${\ast}$}--ideals
$I\neq I'$ of $D$ are such that $(IT)^\ast \neq (I'T)^\ast$.     The
second part  of the statement    follows immediately from the fact
that,  if we set $\ast := \dot{\star}$, then \d{$\ast$} = $\star$
  (Proposition \ref{pr:1.2}(2)).
\end{proof}

\begin{rem}
The converse of  (2) in    Lemma \ref{lemma:noeth} does not hold in
general.  For instance, take $D  \subset    T$, where $D$ is a
Noetherian domain and $T$ is a non-Noetherian overring of $D$.     Let
 $\ast := d_{T}$ and    $\star :=  \star_{\{T\}}$.   Note
that \d{$\ast$} = $\star$.     Then, $D$ is $\star$--Noetherian, by
(1a)    of Lemma \ref{lemma:noeth},    but $D^\star= T^\ast =   T$
is not  $\ast$--Noetherian   (or, equivalently,
$\dot{\star}^{T}$--Noetherian),   because  $\ast =d_T \ \left(
=\dot{\star}^{T} = \dot{\star} \right)$   and $T$ is not Noetherian. \\
However,  if
$\star= \tilde{\star}$, the last statement of (2) in Lemma \ref{lemma:noeth} can be reversed,
as we will see in Proposition \ref{prop:3.5}.
\end{rem}

\begin{lemma} \label{lemma:3.1}
Let $D$ be an integral domain and let $\star$ be a semistar
operation on $D$. Then, $D$ is a $\star$--Noetherian domain if and
only if, for each nonzero ideal $I$ of $D$,     there exists a
finitely generated ideal    $J\subseteq I$  of $D$ such that
$I^\star=J^\star$. Therefore, $D$ is a $\star$--Noetherian domain if and
only if, for each $E\in \boldsymbol{F}(D)$, there
    exists $F \in \boldsymbol{f}(D)$, such that $F \subseteq E$ and
    $F^\star = E^\star$. In particular, if $\star$ is a \sl star
    operation \it on $D$ and if $D$ is a $\star$--Noetherian
    then  $\star$ is a star operation of finite type on $D$.
\end{lemma}
\begin{proof} For the ``only if'' part, let $x_1\in I$, $x_1\neq 0$, and
set $I_1  :=   x_1D$.  If $I^\star=I_1^\star$ we are done.
Otherwise, it is easy to see that    $I \not \subseteq {I_1}^\star
\cap D$.  Let $x_2\in I\smallsetminus \left(I_1^\star \cap D\right)$ and set
$I_2  :=   (x_1, x_2)D$.  By iterating this process, we construct a chain
$\{I_n^\star \cap D \}_{n\ge 1}$ of quasi--${\star}$--ideals of $D$.  By
assumption this chain    must stop, i.e., for some $ k     \geq 1$,
$I^\star _{ k    }\cap D =I^\star _{ k    +1}\cap D$,  and so
   $I^\star_{ k    } = (I^\star_{ k    } \cap D)^\star
=I^\star$.  So, we conclude  by taking $J:= I_{k}$.     Conversely,
let $\{I_n\}_{n\ge 1}$ be a chain of quasi--$\star$--ideals in $D$ and
set $I:=\bigcup_{n\ge 1}I_{ n  }$.  Let $J \subseteq I$ be a finitely
generated ideal of $D$ such that $J^\star=I^\star$, so there exists $k
\ge 1$ such that $J \subseteq I_k$ and $J^\star=I_k^\star=I^\star$.
This implies that the chain of quasi--$\star$--ideals $\{I_n\}_{n\ge 1}$
stops (in fact,  $I_n=I_{k}=I^\star \cap D$, for each $n\geq
k$).
\end{proof}

\begin{prop}\label{prop:3.5}
    Let $D$ be an integral domain and let $\star$
be a semistar operation on $D$.
\begin{enumerate}
      \item[(1)] Assume that $\star$ is stable. Then  $D$ is
${\star}$--Noetherian if and only if $D^{{\star}}$ is
${\dot{{\star}}}$--Noetherian.
\item[(2)] $D$ is
$\tilde{\star}$--Noetherian if and only if $D^{\tilde{\star}}$ is
${\dot{ \tilde{\star}}}$--Noetherian.
\end{enumerate}
\end{prop}
\begin{proof}   (1)    The ``if" part follows from Lemma
\ref{lemma:noeth}(2) and  Proposition \ref{pr:1.2}(2)  (without using the hypothesis of
stability).    Conversely,
    let $I$ be a nonzero ideal of
$D^{{\star}}$ and set $J:=I\cap D$.   Then,
$J^\star = (I \cap D)^\star = I^\star \cap D^\star = I^\star$. Therefore,  by Lemma
\ref{lemma:3.1} (applied to $D$),
we can find $F \in \boldsymbol{f}(D)$ such that $F \subseteq J$ and
$F^\star = J^{\star}$. Hence, $(FD^\star)^{\dot{{\star}}}= F^\star = J^\star =
I^\star = I^{\dot{{\star}}}$. The conclusion follows from Lemma
\ref{lemma:3.1} (applied to $D^{{{\star}}}$, since $FD^\star
\subseteq I$ and $FD^\star \in \boldsymbol{f}(D^\star)$). \\
(2) is a straightforward consequence of (1).
\end{proof}

\begin{prop} \label{prop:noeth}
Let $D$ be an integral domain and $\star$ a semistar operation on
$D$. Then, $D$ is $\star$--Noetherian if and only if $D$ is
$\star_{_{\!f}}$--Noetherian.
\end{prop}
\begin{proof} The ``if'' part follows from Lemma
\ref{lemma:noeth}(1), since $\star_{_{\!f}} \le \star$.     The converse
follows immediately from Lemma \ref{lemma:3.1}.
\end{proof}

\begin{rem} \label{rk:3.6}  Let $D$ be an integral domain
	and $\star$ a semistar operation on $D$. 
	
	\bf (1) \rm  Let $E\in \overline{\boldsymbol{F}}(D)$,
	we say that $E$ is \it $\star$--finite \rm if there exists $F\in \boldsymbol{f}(D)$
	such that $E^\star=F^\star$. 
	From Lemma \ref{lemma:3.1} it follows that if $D$ is a $\star$--Noetherian domain, 
	then each nonzero fractional ideal is $\star$--finite. 
	The converse does not hold in general \cite[Example 18]{GJS99}.
	However, when $\star= \star_{_{\!f}}$, 
	 $E\in \overline{\boldsymbol{F}}(D)$ 
 is  $\star$--finite if and only if there exists $F\in \boldsymbol{f}(D)$
	such that $E^\star=F^\star$, with $F \subseteq E$ \cite[Lemma 
	2.3]{FP} (note that the star operation case was investigated in 
	\cite{AA}).
	From the previous considerations, from Lemma \ref{lemma:3.1} and from 
	Proposition \ref{prop:noeth},
	we deduce easily that \sl
 $D$ is a $\star$-Noetherian domain
 if and only if 
every nonzero fractional ideal of $D$ is  $\star_{_{\!f}}$-finite. \rm

%
%
%
%
   
\bf (2) \rm Note that:\\
\centerline{$\tilde{\star}\mbox{--Noetherian}\; \Rightarrow \;
\star\mbox{--Noetherian},$}\\
because $\tilde{\star} \le \star$  (Lemma \ref{lemma:noeth}(1)).    The
converse is not true in general.  Indeed, if $\star  :=    v$, then
$\star_{_{\!f}} = t$ and $\tilde{\star} = w$ and  we know that
$v$--Noetherian  (= $t$--Noetherian)   is Mori and  that
   $w$--Noetherian is strong Mori  \cite[Section 4]{WMc97}.     Since
it is possible to give examples of Mori domains that are not strong Mori
 \cite[Corollary 1.11]{WMc99},    we  deduce    that
$\star$--Noetherian does not imply $\tilde{\star}$--Noetherian.
\end{rem}

  In the next result, we provide a sufficient condition for the transfer
of the semistar Noetherianity to overrings.

\begin{prop}\label{prop:3.8}
Let $D$ be an integral domain and let $T$ be an overring of $D$.
Let $\star$ be a semistar operation on $D$ and $\star^\prime$ a
semistar operation on $T$. Assume that $T$ is
$(\star,\star^\prime)$--flat over $D$. If $D$ is
$\tilde{\star}$--Noetherian, then $T$ is
$\widetilde{\star^\prime}$--Noetherian.
\end{prop}

\begin{proof}
Let $A$ be a nonzero ideal of $T$. Let $N \in
\mathcal{M}(\widetilde{\star^\prime})=\mathcal{M}({\star^\prime}_f)$
(Proposition \ref{prop:nagata}(5)).
From the $(\star,\star^\prime)$--flatness, it follows that
$T_N=D_{N \cap D}$. Then, $A^{\widetilde{\star^\prime}}=  \cap \{AT_N
\mid  N \in \mathcal{M}({\star^\prime}_f)\} = \cap \{AD_{N \cap
D} \mid  N \in \mathcal{M}({\star^\prime}_f) \}$.  Now, $N \cap D$
is a prime of $D$ such that $(N \cap D)^{\tilde{\star}} \neq
D^{\tilde{\star}}$ (by \cite[Proposition 3.2]{EF}, since $T$ is
$(\star,\star^\prime)$--linked to $D$, by definition of
$(\star,\star^\prime)$--flatness). Hence, $N \cap D$ is a
quasi--$\tilde{\star}$--ideal. Consider the ideal $A \cap D$ of $D$. Since
$D$ is  $\widetilde{\star}$--Noetherian, 
it follows by Lemma \ref{lemma:3.1}
that there exists a finitely generated ideal $C$ of $D$, such that
$C \subseteq A \cap D$ and $C^{\tilde{\star}} = (A\cap
D)^{\tilde{\star}}$. Then, $AT_N=AD_{N \cap D} = (A \cap D)
D_{N\cap D} = (A \cap D)^{\tilde{\star}} D_{N \cap D} =
C^{\tilde{\star}} D_{N \cap D} = C D_{N \cap D} = (CT)T_N$.  Thus, 
$A^{\tilde{\star'}} = (CT)^{\tilde{\star'}}$,  with $CT$ finitely
generated ideal of $T$, such that $CT \subseteq A$. Hence, $T$ is
$\widetilde{\star^\prime}$--Noetherian.
\end{proof}

 Let $D$ be an integral domain and $\star$ a semistar operation on
$D$.     We say that $D$ has the  \it $\star$--finite character
property \rm (for short, \it $\star$--FC property\rm )    if each nonzero
element $x$ of $D$ belongs to only finitely many quasi--$\star$--maximal
ideals of $D$.   Note that the $\star_{_{\!f}}$--FC property
coincides with the $\tilde{\star}$--FC property, because
$\mathcal{M}(\star_{_{\!f}}) = \mathcal{M}(\tilde{\star})$
  (Proposition\ref{prop:nagata}(5)).

\begin{prop} \label{prop:GP1.8}
     Let $D$ be an integral domain and $\star$ a semistar operation on
$D$.     If $D$ is  $\tilde{\star}$--Noetherian,    then $D_M$ is
Noetherian, for each $M \in \mathcal{M}(\star_{_{\!f}})$.  Moreover, if
$D$ has the $\star_{_{\!f}}$--FC property, then the converse
holds.
\end{prop}
\begin{proof} Let $M \in \mathcal{M}(\star_{_{\!f}})$, $A$ an ideal of
$D_M$ and $I  :=   A \cap D$.  Since $D$ is $\tilde{\star}$--Noetherian,
there exists a finitely generated ideal $J \subseteq I$ of $D$ with $J
^{\tilde{\star}} = I^{\tilde{\star}}$ (Lemma \ref{lemma:3.1}).  Then,
$A=ID_M=I^{\tilde{\star}} D_M=J^{\tilde{\star}} D_M = JD_M$  (we used
twice the fact that $\tilde{\star}$ is spectral, defined by
$\mathcal{M}(\star_{_{\!f}})$).     Then $A$ is a finitely generated
 ideal of $D_{M}$    and  so    $D_M$ is Noetherian.  For the
converse, assume that the $\star_{_{\!f}}$--FC property holds  on $D$.
Let $I$ be a nonzero ideal of $D$ and let $0 \ne x \in I$.  Let $M_1,
M_2,   \ldots, M_n \in \mathcal{M}(\star_{_{\!f}})$ be the
quasi--$\star_{_{\!f}}$--maximal ideals containing $x$.  Since $D_{M_i}$
is Noetherian for each $i=1,2, \ldots, n$,  then     $ID_{M_i} = J_i
D_{M_i}$, for some finitely generated ideal $J_i \subseteq I$  of
$D$.   The ideal $B :=   xD + J_1 + J_2 + \ldots + J_n$ of $D$ is finitely
generated and contained in $I$.  It is clear that, for each  $i=1,2,
\ldots, n$,     $ID_{M_i}= BD_{M_i}$.  Moreover, if $M \in
\mathcal{M}(\star_{_{\!f}})$ and $M \neq M_i$,  for each    $i=1, 2, \ldots,n$,
then $x \not \in M$  and this fact    implies $ID_M=BD_M=D_M$.
Then, $I^{\tilde{\star}}=  \bigcap \{ ID_M \mid M \in
\mathcal{M}(\star_{_{\!f}}) \}=\bigcap \{ BD_M \mid M \in
\mathcal{M}(\star_{_{\!f}}) \}=    B^{\tilde{\star}}$.  Thus, by Lemma
\ref{lemma:3.1}, $D$ is $\tilde{\star}$--Noetherian.
\end{proof}

\begin{rem}   \bf (1) \rm Note that Proposition \ref{prop:GP1.8}, 
in case of star operations, can be deduced from \cite[Proposition 
4.6]{HK:2000}, proven in the context of weak ideal systems on 
commutative monoids.

	\bf (2) \rm  Note that strong Mori domains (that is,  $w$--Noetherian domains,
 where $w :=\tilde{v}$) or, more generally, Mori domains    satisfy
 always the $t$--FC property\ (= $w$--FC property, since $\mathcal{M}(w)
 =\mathcal{M}(t)$, for every integral domain)   by  \cite[Proposition
 2.2(b)]{BG87}.     But it is not true in general that  the
 $\tilde{\star}$--Noetherian domains satisfy the $\star_{_{\!f}}$--FC
 property (take, for instance, $ D:=\mathbb{Z}[X]$, $\star :=
 d$,   and observe that $X$ is contained in infinitely many maximal
 ideals of $\mathbb{Z}[X]$). \\
  Note that, from Proposition \ref{prop:GP1.8}  and from the previous
  considerations, we  obtain in particular
 that \sl  an integral domain $D$ is strong Mori if and only if $D_M$ is
 Noetherian, for each $M \in \mathcal{M}(t)$,
 and  $D$ has the $w$--FC property \rm  (cf. also \cite[Theorem
 1.9]{WMc99}).
\end{rem}

%
%

\section{Semistar Dedekind domains}
Let $D$ be an integral domain and $\star$ a semistar operation on
$D$.  We recall  from    Section 1 (or \cite[Section 2]{FP})     that a nonzero
fractional ideal $F\  (\in \boldsymbol{F}(D))$   of $D$ is
\emph{$\star$--invertible} if $(FF^{-1})^\star=D^\star$ and  $E\ \in
\boldsymbol{\overline{F}}(D)$ is \emph{quasi--$\star$--invertible} if
$(E(D^\star:E))^\star=D^\star$ (note that,  the last property
implies that   $E\in
\boldsymbol{F}(D^\star)$).     It is clear that a $\star$--invertible
ideal is quasi--$\star$--invertible.  The converse is not true in
general  \cite[Example 2.9 and Proposition 2.16]{FP}   but, if
$\star$ is stable (e.g., for   $\star=   \tilde{\star}$),   a
finitely generated ideal is $\star$--invertible if and only if it is
quasi--$\star$--invertible \cite[Corollary 2.17(2)]{FP}.

\begin{prop}  \label{prop:3.2} Let $D$ be an integral domain and
$\star$ a semistar operation on $D$. The following are equivalent:
\begin{enumerate} \item[(1)] $D$ is a $\star$--Noetherian domain and a
P$\star$MD;   \item[(1$_{_{\!f}}$)] $D$ is a
$\star_{_{\!f}}$--Noetherian domain and a P$\star_{_{\!f}}$MD;
\item[(2)]
$\boldsymbol{F}^{\tilde{\star}}(D)
:=   \{F^{\tilde{\star}}\mid \,  F\in \boldsymbol{F}(D)\}$ is a group under the
multiplication  ``$\times$'', defined by    $F^{\tilde{\star}}\times
G^{\tilde{\star}} :=
(F^{\tilde{\star}}G^{\tilde{\star}})^{\tilde{\star}}=(FG)^{\tilde{\star}}$,
 for all  $F, G \in \boldsymbol{F}(D)$;   \item[(3)] Each nonzero
fractional ideal  of $D$    is quasi--$\tilde{\star}$--invertible;
\item[(4)] Each nonzero  (integral)    ideal  of $D$    is
quasi--$\tilde{\star}$--invertible.
\end{enumerate}
\end{prop}
\begin{proof}  (1) $\Leftrightarrow$ (1$_{_{\!f}}$) is
obvious   (Proposition  \ref{prop:noeth} and Proposition
\ref{pr:1.5} ((i)$\Leftrightarrow$(vi))).        \\
    (1) $\Rightarrow$ (2).  One can easily check that
$\boldsymbol{F}^{\tilde{\star}}(D)$ is a monoid, with
$D^{\tilde{\star}}$ as the identity element  (with respect to
``$\times$'').     We next show that each element of
$\boldsymbol{F}^{\tilde{\star}}(D)$ is invertible for the monoid
structure.  Let $F \in \boldsymbol{F}(D)$, then there exists $0 \neq d
\in D$ such that $I  :=   dF \subseteq D$.  Write
$I^{\star_{_{\!f}}}=J^{\star_{_{\!f}}}$, where $J \subseteq I$ is a
finitely generated ideal of $D$ (Lemma \ref{lemma:3.1} and Proposition
\ref{prop:noeth}).  Since $D$ is a P$\star$MD, then
$\star_{_{\!f}}=\tilde{\star}$
  (Proposition  \ref{pr:1.5}).
So,
$I^{\tilde{\star}}=J^{\tilde{\star}}$.  We have
$(JJ^{-1})^{\tilde{\star}}=D^{\tilde{\star}}$, since $D$ is a
P$\tilde{\star}$MD
  (Proposition  \ref{pr:1.5}).
Then,
$D^{\tilde{\star}}=(J^{\tilde{\star}}J^{-1})^{\tilde{\star}}=
(IJ^{-1})^{\tilde{\star}}=(dFJ^{-1})^{\tilde{\star}}=
\left(F^{\tilde{\star}}(dJ^{-1})^{\tilde{\star}}\right)^{\tilde{\star}}$.
Thus $F^{\tilde{\star}}$ is invertible  in
$(\boldsymbol{F}^{\tilde{\star}}(D), \times)$.     \\
(2) $\Rightarrow$ (3).  Let $F \in \boldsymbol{F}(D)$.   By
assumption,   there exists $G\in \boldsymbol{F}(D)$ such that
$(FG)^{\tilde{\star}} = D^{\tilde{\star}}$.  We have $FG \subseteq
D^{\tilde{\star}}$, so $G \subseteq (D^{\tilde{\star}}:F)$.  Thus
$D^{\tilde{\star}}= (FG)^{\tilde{\star}} \subseteq
(F(D^{\tilde{\star}}:F)) \subseteq D^{\tilde{\star}}$.  Hence
$(F(D^{\tilde{\star}}:F))^{\tilde{\star}}=D^{\tilde{\star}}$, that is,
$F$ is quasi--$\tilde{\star}$--invertible.  \\
(3) $\Rightarrow$ (4) is straightforward. \\
(4) $\Rightarrow$ (1) From the  previous    comments on quasi semistar
invertibility  for nonzero finitely generated ideals in the stable case, it
is clear that the assumption implies that    $D$ is a
P$\tilde{\star}$MD and hence $D$ is a P$\star$MD
  (Proposition  \ref{pr:1.5}).
To prove that $D$ is a $\star$--Noetherian domain, since
$\tilde{\star}=\star_{_{\!f}}$
  (Proposition  \ref{pr:1.5}),
it is enough to show,  by using
Proposition \ref{prop:noeth},    that $D$ is
$\tilde{\star}$--Noetherian.  Let $I$ be a nonzero ideal of $D$,
then, by assumption,
$(I(D^{\tilde{\star}}:I))^{\tilde{\star}}=D^{\tilde{\star}}$.   By
\cite[Lemma 2.3 and Proposition 2.15]{FP} applied to ${\tilde{\star}}$,
   there exists a nonzero finitely generated ideal $J$ of $D$ such that
 $J \subseteq I$ and $ J^{\tilde{\star}} = I^{\tilde{\star}}$.
From Lemma \ref{lemma:3.1}, we deduce that $D$ is
 $\tilde{\star}$--Noetherian.
\end{proof}

An integral domain $D$   with a semistar operation
$\star$ satisfying      the equivalent conditions (1)--(4) of Proposition
\ref{prop:3.2} is called a \emph{$\star$--Dedekind domain}
(\emph{$\star$--DD}, for short).   Note that, by definition, the notions of
$\star$--DD and $\star_{_{\!f}}$--DD coincide.

\begin{rem} \label{rem:krull}
(1)  By Proposition
\ref{prop:3.2}(1),    if $\star=d$ we obtain  that a $d$--DD coincides
with    a classical Dedekind domain  \cite[Theorem 37.1]{Gil92}; \ if $\star=v$, we have  that a
$v$--DD coincides with a Krull domain (since a Mori P$v$MD is a Krull
domain \cite[Theorem 3.2 ((1) $\Leftrightarrow$(3))]{K89}; 
note that a Mori domain verifies the $t$--FC property by
\cite[Proposition 2.2(b)]{BG87}).   More generally, if $\star$ is 
a star operation, then $D$ is a $\star$--DD if and only if $D$ is  
$\star$--Dedekind in the sense of \cite[Chapter 23]{HK98}. 

 (2)    If $D$ is $\star$--DD then   $D$ is      $\star$--ADD  (for a converse, see
the following Theorem \ref{thm:3.3}).     Indeed,  since a
$\star$--DD is a P$\star$MD and so    $\tilde{\star}=\star_{_{\!f}}$
  (Proposition  \ref{pr:1.5}).
 This
equality implies also that    $D$ is $\tilde{\star}$--Noetherian
(Proposition \ref{prop:noeth} and Proposition \ref{prop:3.2}(1)).
Therefore    $D_M$ is Noetherian (by Proposition \ref{prop:GP1.8}) and,
hence,  we conclude that    $D_M$ is a DVR,  for each $M
\in \mathcal{M}(\star_{_{\!f}})$.

\end{rem}

\begin{cor} \label{cor:tilde} Let $D$ be an integral domain and $\star$ a semistar
operation on $D$. Then $D$ is a $\star$--DD if and only if $D$ is a
$\tilde{\star}$--DD.
\end{cor}
\begin{proof} It follows from Proposition  \ref{prop:3.2}(4)    and
from   the fact that $\tilde{\tilde{\star}}=\tilde{\star}$,
 since $\mathcal{M}(\widetilde{\star})= \mathcal{M}(\star_{_{\!f}})$
(cf. also \cite[page 182]{FH2000}).
\end{proof}

\begin{thm} \label{lemma:krull}  Let $D$ be an integral domain.
\begin{enumerate}
\item[(1)] Let $\star \leq \star'$ be two semistar operations on
$D$. Then:

\centerline{$D$  is a    $\star\mbox{--}DD \; \Rightarrow \;  D$
is a    $\star'\mbox{--}DD$\,.}

In particular:
\begin{enumerate}
 \item[(1a)]    If $D$ is a Dedekind domain, then $D$ is a $\star$--DD,
for any semistar operation $\star$ on $D$.
 \item[(1b)]      Assume that    $\star$ is a (semi)star operation  on $D$.
  Then a    $\star$--DD is a Krull domain.
\end{enumerate}
\item[(2)]  Let $T$ be an overring of $D$.  Let $\star$ be a semistar
operation on $D$ and $\star'$ a semistar operation on $T$.  Assume that
$T$ is a $(\star, \star')$--linked overring of $D$.  If $D$ is a
$\star$--DD, then $T$ is a $\star'$--DD. In particular, If $D$ is a
$\star$--DD, then $D^\star$ is a $\dot{\star}$--DD.
\end{enumerate}
\end{thm}
\begin{proof}
(1) follows from \cite[p. 30]{FJS03} and Lemma \ref{lemma:noeth}(1).
 (1a) and (1b) are consequence of (1), Remark \ref{rem:krull}(1) and
 of the fact that $d \leq \star$, for each semistar operation
$\star$, and if $\star$ is a (semi)star operation, then $\star \leq
v$.   \\
(2)
  Note that if $T$ is a $(\star, \star')$--linked overring of $D$
and if $D$ is a P$\star$MD, then  $T$ is a $(\star,
\star')$--flat over $D$ \cite[Theorem 5.7 ((i)$\Rightarrow$(ii))]{EF}.
By Proposition \ref{prop:3.2}(1) and Corollary  \ref{cor:tilde}, we
know that $D$  is $\tilde{\star}$--Noetherian and a
P${\star}$MD (or, equivalently, a P$\tilde{\star}$MD). Hence, $T$ is $\widetilde{\star'}$--Noetherian
(Proposition \ref{prop:3.8}) and $T$ is a P${\star'}$MD (or, equivalently,
a P$\tilde{\star'}$MD) by \cite[Corollary 5.4]{EF}.  The  first
statement 
follows from Proposition \ref{prop:3.2}(1) and Corollary
\ref{cor:tilde}.  The last statement is a consequence of
\cite[Lemma 3.1(e)]{EF}. 
\end{proof}

\begin{prop} \label{prop:krull} Let $D$ be an integral domain and
$\star$ a \sl (semi)star operation \it on $D$.  Then the following are
equivalent:
\begin{enumerate}
\item[(1)] $D$ is a $\star$--DD
\item[(2)] $D$ is a Krull domain and $\star_{_{\!f}}=t$
\end{enumerate}
\end{prop}
\begin{proof}  (1) $\Rightarrow$ (2).     By Theorem \ref{lemma:krull}(1b),
if $D$ is a $\star$--DD, then $D$ is a Krull domain, in this case,
$\star_{_{\!f}}=t$ \cite[Proposition 3.4]{FJS03}.\\
  (2) $\Rightarrow$ (1).  This follows from Remark
  \ref{rem:krull}(1)   and from the fact that   $v$--DD = $t$--DD = $\star_{_{\!f}}$--DD = $\star$--DD.
\end{proof}

Note that  Proposition \ref{prop:krull}   has already been
proven in \cite[Theorem 23.3((a)$\Leftrightarrow$(d))]{HK98}, by using the language of monoids and
ideal systems.

 \begin{rem} Note that if $D$ is $\star$--DD, then by  Theorem
 \ref{lemma:krull}(2) $D^\star$ is $\dot{\star}$--DD, that is   $D^\star$ is a Krull
 domain and
 $(\dot{\star})_f=t_{D^\star}$  (Proposition \ref{prop:krull}).   However, the converse does not hold in
 general as the example in Remark \ref{rem:2.1}(2) shows.  Nevertheless,
  the converse is true when $\star =\tilde{\star}$ (see the following
 Corollary \ref{cor:3.8}) or when the extension $D\subseteq D^\star$ is flat,
 as a consequence of Lemma \ref{lemma:noeth}(2)  and  \cite[Proposition
 3.2]{FJS03}.
 For a more accurate
 discussion on this problem see the following Remark \ref{rk:4.17}.   \par
     \end{rem}

     Next result is a ``Cohen-type'' Theorem for quasi--$\star$--invertible
     ideals.

\begin{lemma} \label{lemma:GP3.4} Let $D$ be an integral domain and $\star$ a semistar
operation of finite type  on $D$.    The following are equivalent:
\begin{enumerate}
\item[(1)] Each
nonzero    quasi--$\star$--prime  of $D$    is a
quasi--$\star$--invertible ideal  of $D$.   \item[(2)] Each nonzero
 quasi--$\star$--ideal of $D$ is a quasi--$\star$--invertible
ideal  of $D$.    \item[(3)] Each nonzero ideal of $D$ is a
quasi--$\star$--invertible ideal  of $D$.
\end{enumerate}
\end{lemma}
\begin{proof}
    (1) $\Rightarrow$ (2).  Let $S$ be the set
    of the quasi--$\star$--ideals of $D$ that are not
    quasi--$\star$--invertible.  Assume that $S \neq \emptyset$.   Since
    $\star = \star _{_{\!f}}$ by assumption, then Zorn's Lemma can be
    applied, thus we deduce that    $S$ has maximal elements.  We next show
    that a maximal element of $S$ is prime.  Let $P$ be a maximal element of
    $S$ and let $r,s \in D$, with $rs \in P$.  Suppose $s \not \in P$.  Let
    $J :=   (P:_DrD)$.  We claim that $J^{\star} \cap D = J$.  Indeed,
    since $(P:_DrD)^{\star} \subseteq (P^{\star}:_{D^{\star}}rD)$, then
    $J^{\star}    \cap D \subseteq (P^{\star}:_{D^{\star}}rD) \cap D =
    (P^\star:_DrD)$.  Moreover, if $x \in (P^{\star}:_DrD)$, then $xr \in
    P^{\star} \cap D=P$, and hence $(P^{\star}:_DrD)  \subseteq (P:_DrD)
       =J$.   Thus $J= J^{\star} \cap D$, i.e. $J$ is a
    quasi--$\star$--ideal  of $D$.     Clearly, $J$ contains properly
    $P$ (since $s \in J\smallsetminus P$).   By the maximality of $P$ in
    $S$,    it follows that $J$ is quasi--$\star$--invertible, that is
    $(J(D^{\star}:J))^{\star}=D^{\star}$.  We notice that $P(D^{\star}:J)\in
    \boldsymbol{\overline{F}}(D)$ is not quasi--$\star$--invertible, since
    $P$ is not quasi--$\star$--invertible  \cite[Lemma 2.10]{FP}.     We
    deduce that $(P(D^{\star}:J))^{\star} \cap D $ is a  proper
    quasi--$\star$--ideal, that is not quasi--$\star$--invertible
    \cite[Remark 2.13(a)]{FP} and, obviously, it    contains $P$.  From the
    maximality of $P$  in $S$,   we have $(P(D^{\star}:J))^{\star} \cap
    D = P$.  Now, $rJ \subseteq P$ implies $(rJ)^{\star} \subseteq
    P^{\star}$.  Then $r \in (rD)^{\star}  =
    (rJ(D^{\star}:J))^{\star} \subseteq (P(D^{\star}:J))^{\star}$.
    Therefore, $r \in (P(D^{\star}:J))^{\star} \cap D = P$ and  so we
    have proven that    $P$ is  a prime ideal of $D$.     \\
 (2) $\Rightarrow$ (3) is a consequence of \cite[Remark 2.13(a)]{FP},
after remarking that, for each nonzero ideal $J$ of $D$, then $J
\subseteq I:= J^\star \cap D$,   where    $I$ is a quasi--$\star$--ideal of
$D$ and $J^\star = I^\star$. \\
    (3) $\Rightarrow$    (2) $\Rightarrow$ (1) are trivial.
\end{proof}

\begin{rem} Note that, in the situation of Lemma \ref{lemma:GP3.4}, the
statement:
\begin{enumerate}
\item[(0)] \it each nonzero quasi--$\star$--maximal ideal of $D$ is a
quasi--$\star$--invertible ideal of $D$, \rm
\end{enumerate}
is, in general, strictly weaker than (1).  Take, for instance, $D$ equal to a
discrete valuation domain of rank $\geq 2$, and $\star = d_{D}$.
\end{rem}

   The next   two   theorems    generalize \cite[Theorem 37.8
((1)$\Leftrightarrow$(4)), Theorem 37.2]{Gil92}.  Similar results are proven in
\cite[Theorem 23.3((a)$\Leftrightarrow$(c), (h))]{HK98}.

\begin{thm} \label{thm:GP3.4} Let $D$ be an integral domain and $\star$ a semistar
operation on $D$. The following are equivalent:
\begin{enumerate}
\item[(1)] $D$ is a $\star$--DD.
\item[(2)] Each  nonzero    quasi--$\tilde{\star}$--prime  ideal of
$D$    is quasi--$\tilde{\star}$--invertible.
\end{enumerate}
\end{thm}
\begin{proof}
  Easy consequence of    Lemma \ref{lemma:GP3.4}
((1)$\Leftrightarrow$(3))   and
Proposition \ref{prop:3.2} (4).
\end{proof}

 From the previous theorem, we deduce the following
characterization of Krull domains   (cf. \cite[Theorem 2.3
((1)$\Leftrightarrow$(3))]{HZ}, \cite[Theorem 3.6
((1)$\Leftrightarrow$(4))]{K89} and
\cite[Theorem 5.4 ((i)$\Leftrightarrow$(vi))]{WMc97}).

\begin{cor} \label{cor:krull}  Let $D$ be an integral domain. The following are equivalent:
\begin{enumerate}
\item[(1)] $D$ is a  Krull domain.
\item[(2)] Each  nonzero $w$--prime ideal of
$D$ is $w$--invertible.
\item[(3)] Each  nonzero $t$--prime ideal of
$D$ is $t$--invertible.
\end{enumerate}
\end{cor}
\begin{proof} (1) $\Leftrightarrow$ (2) is a direct consequence of
Theorem \ref{thm:GP3.4}.\\
(1) $\Rightarrow$  (3)  is a straightforward consequence of  (1) $\Rightarrow$
(2)  and of the fact that, in a Krull domain (which is a particular
P$v$MD), $t= \widetilde{t}=w$
(Proposition \ref{pr:1.5}).  \\
(3) $\Rightarrow$ (2).  Note that, by assumption, and by Lemma \ref{lemma:GP3.4}
((1)$\Leftrightarrow$(3)), every nonzero ideal of $D$ is
$t$--invertible.  Let $Q$ be a nonzero
$w$--prime.  If  $(QQ^{-1})^w \neq D$, then $Q
\subseteq (QQ^{-1})^w \subseteq M$, for
some $M \in  \mathcal{M}(w) =
\mathcal{M}(t) $ (Proposition \ref{prop:nagata}(5)), thus $
(QQ^{-1})^t = ((QQ^{-1})^w)^t \subseteq M^t =M$, which is a
contradiction.
   \end{proof}

\begin{thm} \label{thm:3.3}  Let $D$ be an integral domain and $\star$
a semistar operation on $D$.     The following are equivalent:
\begin{enumerate}
\item[(1)] $D$ is a $\star$--DD.
\item[(2)] $D$ is  a    $\star$--ADD and each nonzero element of $D$ is
contained in only finitely many quasi--$\star_{_{\!f}}$--maximal ideals
(i.e. $D$ has the $\star_{_{\!f}}$--FC property).
\item[(3)] $D$ is  a  $\star$--Noetherian  $\star$--ADD.
\end{enumerate}
\end{thm}
\begin{proof} (1) $\Rightarrow$ (2). Clearly $D$ is a $\star$--ADD, by Remark \ref{rem:krull}(2). Since
each quasi--$\star_{_{\!f}}$--maximal ideal of $D$ is a contraction of a
$\dot{\star}_f$--maximal ideal of $D^\star$ \cite[Lemma
2.3(3)]{FL03},  in order    to show that $D$ has $\star_{_{\!f}}$--FC
property, it is enough to check that $D^\star$ satisfies the
$\dot{\star}_{_{\!f}}$--FC property.   On the other hand,    since (1)
implies that $D^\star$ is a $\dot{\star}$-DD  (Theorem
\ref{lemma:krull}(2)), without loss of generality,    we can assume
that $\star$ is a (semi)star operation on $D$ and $D$ is a $\star$--DD.
By Proposition \ref{prop:krull}, $D$ is a Krull domain and
$\star_{_{\!f}}=t$.  Thus, each nonzero element is contained in only
finitely many $t$--maximal ideals ($=\star_{_{\!f}}$--maximal ideals)
 of $D$.     \\
 (2) $\Rightarrow$ (1).   We  need to show that
$D$ is ${\star_{_{\!f}}}$--DD. First, note that $D$ is a P${\star_{_{\!f}}}$MD
and $D_M$ is Noetherian, for each $M\in
\mathcal{M}({\star_{_{\!f}}})$\
(Proposition \ref{prop:2.2} (1) and (2)).  The conclusion now follows
from Proposition \ref{prop:GP1.8} and Proposition
\ref{prop:3.2}(1),  after recalling that,  in a
P${\star_{_{\!f}}}$MD, $\star_{_{\!f}} =\widetilde{\star}$
(Proposition \ref{pr:1.5}).\\
(1) $\Leftrightarrow$ (3) is a consequence of Proposition
\ref{prop:2.2}(2), Proposition \ref{prop:3.2} and Remark
\ref{rem:krull}(2).
\end{proof}

From the previous theorem, we deduce a restatement of a wellknown characterization of Krull
domains:

\begin{cor} Let $D$ be an integral domain, then the
following are equivalent:
\begin{enumerate}
    \item[(1)] $D$ is a Krull domain.
    \item[(2)] $D$ is a $t$--almost Dedekind domain
    and each nonzero element of $D$ is contained in only finitely many
    $t$--maximal ideals (= $t$--FC property).  \hfill $\Box$
    \end{enumerate}
    \end{cor}

Let $D$ be an integral domain and $\star$ a semistar operation on
$D$.  We recall that the \it   $\star$--integral    closure $D^{[\star]}$
of $D$ \rm (or, the  \it semistar integral closure with respect to the semistar
operation $\star$ of $D$\rm) is the integrally closed overring of
$D^\star$ defined by $D^{[\star]} :=    \{(F^\star:F^\star)  \mid
   F \in \boldsymbol{f}(D)\}$ \cite[Definition 4.1]{FL01}.   We say that
\it $D$ is quasi--$\star$--integrally closed \rm  (respectively,  \it
$\star$--integrally closed \rm) \rm if $D^{\star}=D^{[\star]}$
(respectively, $D=D^{[\star]}$).
 It is clear that:

 --\  $D$ is
quasi--$\star$--integrally closed if and only if $D$ is
quasi--$\star_{_{\!f}}$--integrally closed (respectively,  $D$ is
$\star$--integrally closed if and only if $D$ is
$\star_{_{\!f}}$--integrally closed);

--\   $D$   is    $\star$--integrally closed  if and only if  $D$ is
quasi--$\star$--integrally closed and $\star$ is a
(semi)star operation on $D$.

Note that when $\star =v$, then the overring $D^{[v]}= D^{[t]}$
was studied in \cite{AHZ91} under the name of
\it psedo-integral closure of $D$. \rm

\begin{lemma} \label{lemma:intclos}
Let $D$ be an integral domain and $\star$ a semistar operation on
$D$.
\begin{enumerate}
\item[(1)] If $\star$ is e.a.b., then $D^\star=D^{[\star]}$  (i.e.
$D$ is  quasi--${\star}$--integrally   closed). 
\item[(2)] $D$ is  quasi--$\tilde{\star}$--integrally    closed if and only if
$D^{\tilde{\star}}$ is integrally closed.
\end{enumerate}
\end{lemma}
\begin{proof}
(1)  Note that, in general,   $D^\star \subseteq D^{[\star]}$.  For the converse,
let $F \in \boldsymbol{f}(D)$ and let $x \in (F^\star:F^\star)$.  Then,
$xF^\star \subseteq F^\star$ and $F^\star = F^\star + F^\star (xD)$.
 Therefore we have $ (F (D+xD))^\star =   (F^\star (D+xD))^\star =
 (F^\star + F^\star (xD))^{\star} = F^\star$.   From the fact that $F$ is
finitely generated and that $\star$ is e.a.b., we obtain $ (D+xD)^\star
= D^\star $.  It follows that $x \in D^\star$ and  so    $(F^\star :
F^\star) \subseteq D^\star$.  Hence, $D^\star=D^{[\star]}$.\\
(2)  The ``only if'' part is clear.  For the ``if" part, let $D'$ be
the integral closure of $D$, since $D^{\tilde{\star}}$ is integrally
closed, then $(D')^{\tilde{\star}}\subseteq D^{\tilde{\star}}\subseteq
D^{[\tilde{\star}]}$ hence, by \cite[Example 2.1(c2)]{FJS03},
$(D')^{\tilde{\star}}= D^{\tilde{\star}}= D^{[\tilde{\star}]}$.
Therefore, $D$ is quasi--$\tilde{\star}$--integrally closed.
\end{proof}

\begin{cor} \label{cor:intclos}  Let  $\star$ be a semistar operation
on an integral domain $D$.
If $D$ is a P$\star$MD (in particular, a $\star$--DD) then $D$ is
quasi--$\star$--integrally closed.
\end{cor}
\begin{proof} It follows from Lemma \ref{lemma:intclos}(1) and from the fact
 that,    in a P$\star$MD, $\tilde{\star}=\star_{_{\!f}}$ is an
e.a.b. semistar operation    (Proposition \ref{pr:1.5}
((i)$\Rightarrow$(v), (vi))).
\end{proof}

 The following result shows that a semistar version of the ``Noether's
Axioms'' provides a characterization of the semistar Dedekind domains.

\begin{thm} \label{thm:3.5} Let $D$ be an integral domain and $\star$ a semistar operation on $D$.
The following are equivalent:
\begin{enumerate} \item[(1)] $D$ is a $\star$--DD. \item[(2)] $D$ is
$\tilde{\star}$--Noetherian, $\tilde{\star}\mbox{-}\di(D) = 1$ and $D$
is  quasi--$\tilde{\star}$--integrally closed.    \item[(3)] $D$ is
$\tilde{\star}$--Noetherian, $\tilde{\star}\mbox{-}\di(D) = 1$ and
$D^{\tilde{\star}}$ is integrally closed.
\end{enumerate}
\end{thm}
\begin{proof} The equivalence (2) $\Leftrightarrow$ (3) follows from
Lemma \ref{lemma:intclos}  (2).   \\
(1) $\Rightarrow$ (2). Since $D$ is a ${\star}$--DD,
then $D$ is ${\star}$--ADD
 (Remark \ref{rem:krull}(2)).       Hence $\tilde{\star}$-$\di(D)= 1 $
(Proposition \ref{prop:2.8}).   Moreover,    recall that a ${\star}$--DD is a
$\tilde{\star}$--DD (Corollary \ref{cor:tilde}).
Then $D$ is $\tilde{\star}$--Noetherian and a
P$\tilde{\star}$MD (Proposition
\ref{prop:3.2}),
and so    $D$ is
quasi--${\tilde{\star}}$--integrally closed  by Corollary
\ref{cor:intclos}. \\
(3) $\Rightarrow$ (1)  For each $M \in \mathcal{M}(\star_{_{\!f}})$,
it is wellknown that $D^{\tilde{\star}} \subseteq D_{M}$ and $
{D^{\tilde{\star}}}_{MD_M \cap D^{\tilde{\star}}}= D_M$.  Since
$D^{\tilde{\star}}$ is integrally closed, this implies that $D_M$ is also
integrally closed.  Therefore    $D_M$ is a local, Noetherian (by
Proposition \ref{prop:GP1.8}), integrally closed, one dimensional (by
Lemma \ref{lemma:2.5}) domain, that is, a DVR \cite[Theorem
37.8]{Gil92}.  Hence $D$ is a P$\star$MD. In particular, we have
$\tilde{\star}=\star_{_{\!f}}$
  (Proposition \ref{pr:1.5}),
thus    $D$ is
$\star_{_{\!f}}$--Noetherian, by the assumption,  and so $D$ is
$\star$--Noetherian (Proposition \ref{prop:noeth}).   We conclude that
$D$ is a $\star$--DD.
\end{proof}

 By taking $\star=v$ in Theorem \ref{thm:3.5}, we obtain the following
characterization of Krull domains:

\begin{cor} \label{cor:4.16} Let $D$ be an integral domain.  The following are
equivalent:
\begin{enumerate}
    \item[(1)] $D$ is a Krull domain.
    \item[(2)]  $D$ is
 a strong Mori domain, $w\mbox{-}\di(D)=1$ and $D =D^{[w]}$.
 \item[(3)] $D$ is a
 strong Mori domain, $w\mbox{-}\di(D)=1$ and $D$ is integrally closed.
 \item[(4)]$D$
 is a strong Mori domain, $t\mbox{-}\di(D)=1$ and $D$ is integrally
 closed.\end{enumerate}
\end{cor}
\begin{proof}The only part which needs a justification is the statement on
$t$-dimension and $w$-dimension (in the equivalence (3) $\Leftrightarrow$
(4)).  This  follows from the fact that,   in every integral
domain,  $w \leq t$ and 
${\mathcal{M}}(t)={{\mathcal{M}}}(w)$. \end{proof}

\begin{rem}
Note that, if $D$ is a $\star$--DD,  then we know that $\tilde{\star} =
\star_{_{\!f}}$, and so    $D$   satisfies    the properties:
\begin{enumerate}
    \item[($2_{_{\!f}}$)] \sl $D$ is $\star_{_{\!f}}$--Noetherian, $\star_{_{\!f}}\mbox{-}\di(D)=1$ and
$D$ is  quasi--$\star_{_{\!f}}$--integrally closed;
\item[($3_{_{\!f}}$)] \sl  $D$ is
$\star_{_{\!f}}$-Noetherian, $\star_{_{\!f}}\mbox{-}\di(D)=1$ and
$D^{\star_{_{\!f}}}\ (=D^\star)$ is integrally closed \rm
\end{enumerate}
obtained from (2) and (3) of Theorem \ref{thm:3.5}, replacing
$\tilde{\star}$ with $\star_{_{\!f}}$.    But, conversely, if $D$
satisfies  either ($2_{_{\!f}}$) or ($3_{_{\!f}}$)    then $D$ is not
necessarily a $\star$--DD. Indeed, let $D, T$ and $\star$ be as in the
example of Remark \ref{rem:2.1}(2).  Then we have already observed that
$\star = \star_{_{\!f}}$ and $\tilde{\star}=d_D$.  Moreover, $D$ is not
a $\star$--DD (because it is not a $\star$--ADD), but
$D^{\star_{_{\!f}}}=T  =D^{[\star_{_{\!f}}]}$ is integrally closed
(since $T$ is a DVR), $\star_{_{\!f}}$-$\di(D)=1$
(since $\mathcal{M}(\star_{_{\!f}})=\{M\}$   and
$\star_{_{\!f}}$-$\di(D)\leq \di(D) =1$)  and $D$ is
$\star_{_{\!f}}$--Noetherian  (Lemma \ref{lemma:3.1}, since $T$
is Noetherian).

Note that
 ($3_{_{\!f}}$)
does not imply   that $D$ is a $\star$--DD, even if $\star$ is a \sl
(semi)star \rm operation on $D$.  Take $T$ and $D$ as in the example
described in Remark \ref{rem:2.1}(2) and, moreover, assume that $k$ is
algebraically closed in $K$.  It is wellknown that, in this situation,
$D$ is integrally closed.  Let $\star :=v$ on $D$.  It is easy to see
that $\mathcal{M}(v) =\mathcal{M}(t)=\{M\}$, thus $w =d$ is the identity
(semi)star operation on $D$ (hence, $D^{[w]} = D^{[d]}=D$) and
$t\mbox{-}\di(D) =1 \ (= v\mbox{-}\di(D) =w\mbox{-}\di(D) =\di(D))$.
Moreover, it is known that $D$ is a Mori domain \cite[Theorem 4.18]{GH97} and thus
$D$ is a $t$--Noetherian domain.  However, $D$ is not a Krull domain,
since $D$ is not completely integrally closed (being $T$ the complete
integral closure of $D$).  Note that, in this situation, $D$ is even not a
strong Mori domain (by Corollary \ref{cor:4.16}).  \\
 Note also that, in the previous example, $D \subsetneq D^{[t]}$ (i.e.
$D$ is not $t$--integrally closed, hence does not satisfies condition
($2_{_{\!f}}$) for $\star =v$), since  $D^{[t]}= T$ by \cite[Theorem
1.8(ii)]{AHZ91}.

On the other hand,  if $\star$ is a \sl
(semi)star \rm operation on $D$, then we know that $D$ is a
$\star$--DD if and only if  $D$ is a
$v$--DD (i.e. a Krull domain) and $\star_{_{\!f}} =t$ (Proposition
\ref{prop:krull}). It is interesting to observe that, for $\star =v$,  condition (1) of Theorem
\ref{thm:3.5} is equivalent to ($2_{_{\!f}}$). More precisely we have
the following variation of the equivalence (1) $\Leftrightarrow$ (4) of
Corollary \ref{cor:4.16}:\\
\sl $D$ is a Krull domain if and only if $D$ is $t$--Noetherian, $t\mbox{-}\di(D)=1$ and
$D$ is  $t$--integrally closed (i.e. $D = D^{[t]}$). \rm \\
As a matter of fact, let $F \in \boldsymbol{f}(D)$, then  $D =
D^{[t]}=D^{[v]}$ implies that $D =(F^v :F^v) =(F^{-1}
:F^{-1})=(FF^{-1})^{-1}$ and so $(FF^{-1})^{v}=D$.  Moreover, since
$t$--Noetherian is equivalent to $v$--Noetherian (Proposition \ref{prop:noeth}) and $v$--Noetherian
implies that $v =t$ (Lemma \ref{lemma:3.1}), then $(FF^{-1})^{t}=D$. Thus $D$ is a
P$v$MD and so $D$ is a $v$--DD (Proposition \ref{prop:3.2}).

Finally, from the previous considerations we deduce that \sl
$D$ is a
$\star$--DD if and only if

\rm ($\overline{2}_{_{\!f}}$) \sl  $D$ is $\star_{_{\!f}}$--Noetherian,
$\star_{_{\!f}}\mbox{-}\di(D)=1$,
$D$ is  quasi--$\star_{_{\!f}}$--integrally closed and
$\star_{_{\!f}} =t$. \rm

We conclude with a question: is there an example
of an integral (Krull) domain $D$, equipped with a (semi)star
operation $\star$, such that condition  (${2}_{_{\!f}}$)
holds but ($\overline{2}_{_{\!f}}$)  does not? Note that if such an
example exists then necessarily $d \lneq \star_{_{\!f}}\ ( \lneq t )$
\cite[Theorem 37.8 ((1)$\Leftrightarrow$(2))]{Gil92}.

\end{rem}

 Next result generalizes \cite[Proposition 38.7]{Gil92}.

\begin{thm} \label{prop:3.7}
Let $D$ be an integral domain and
$\star$ a semistar operation on $D$. The following are equivalent:
\begin{enumerate}
\item[(1)] $D$ is a $\star$--DD.
\item[(2)] $\Na(D,\star) \ (= \Kr(D, \star))$    is a PID.
 \item[(3)]
$\Na(D,\star) \ (= \Kr(D, \star))$ is a Dedekind domain.
\end{enumerate}
\end{thm}
\begin{proof}
(1) $\Rightarrow$ (2).   Since $D$ is a P$\star$MD, then $\Na(D,\star)=\Kr(D,\star)$
is a B\'ezout domain   (Proposition 1.5 ((i)$\Rightarrow$(iv)) and Proposition 1.4(1)).
Now, let   $\boldsymbol{\mathfrak{I}}$    be a nonzero ideal of $\Na(D, \star)$
and set $I: =\boldsymbol{\mathfrak{I}}\cap D$. We claim that $\boldsymbol{\mathfrak{I}}=
I\!\Na(D, \star)$.
The inclusion $I\!\Na(D, \star)\subseteq \boldsymbol{\mathfrak{I}}$ is clear.
For the opposite inclusion, since $\boldsymbol{\mathfrak{I}}=
(\boldsymbol{\mathfrak{I}}\cap D[X])\Na(D, \star)$,
it is enough to show that $\boldsymbol{\mathfrak{I}}\cap D[X]\subseteq
I\!\Na(D, \star)$. Let $f \in \boldsymbol{\mathfrak{I}}\cap D[X]$,   then
$f\! \Na(D,\star)= f\!\Kr(D,\star)=
{\boldsymbol{c}}(f)\!\Kr(D,\star)={\boldsymbol{c}}(f)\!\Na(D, \star)$
(where the second equality holds by
  Proposition \ref{prop:Kr}(6)).
  Hence ${\boldsymbol{c}}(f)\subseteq f\! \Na(D,\star)\cap D \subseteq
  \boldsymbol{\mathfrak{I}}\cap D=I$.
  Therefore we conclude that $f\in {\boldsymbol{c}}(f)\!\Na(D, \star) \subseteq
  I\!\Na(D,\star)$,
  which proves our claim.
  Now, since $D$ is a $\tilde{\star}$--Noetherian domain (as $D$ is a
  ${\star}$--DD, cf. Corollary \ref{cor:tilde} and Proposition
  \ref{prop:3.2}),
  then  $I^{\tilde{\star}}=
   F^{\tilde{\star}}$ for some  $F \in \boldsymbol{f}(D)$, with
$F\subseteq  I$ (Lemma \ref{lemma:3.1}).  Since $E^{\tilde{\star}}
= E\!\Na(D,\star) \cap K$, for each $E \in {\overline{F}}(D)$(Proposition
\ref{prop:nagata}(4)),
then  we have
$\boldsymbol{\mathfrak{I}}=I\!\Na(D,\star)=I^{\tilde{\star}}\!\Na(D,\star)=
F^{\tilde{\star}}\!\Na(D,\star)= F\!\Na(D,\star)$.
Hence we conclude that $\boldsymbol{\mathfrak{I}}$ is a principal ideal
in $\Na(D,\star)$,
because, as we have already remarked, $\Na(D,\star)$ is a B\'ezout
domain.  \\
 (2) $\Rightarrow$ (3) is trivial. \\
   (3) $\Rightarrow$ (1).
Assume that    $\Na(D, \star)$ is a Dedekind domain then, obviously,
$\Na(D,\star) = \Kr(D, \star)$
(Proposition \ref{pr:1.5} ((i) $\Rightarrow$ (iv)))
and $\Na(D,
\star)$ is an ADD, and hence $D$ is a $\star$--ADD (Theorem
\ref{thm:2.3}).   In order to apply Theorem \ref{thm:3.3}, it remains to show that $D$ has the
$\star_{_{\!f}}$--FC property.   Let $0\neq x\in D$.  Since $\Max(\Na(D,
\star))=\{M\!\Na(D, \star)\mid\, M\in \mathcal{M}(\star_{_{\!f}})\}$
  (Proposition \ref{prop:nagata}(2))
and $\Na(D, \star)$ is  a Dedekind
domain, then    there are only finitely many maximal ideals $M\!\Na(D,
\star)$ containing $x$.  Furthermore, $M\!\Na(D, \star)\cap D=M$, for
each $M\in \mathcal{M}(\star_{_{\!f}})   =\mathcal{M}(\widetilde{\star})$
(Proposition \ref{prop:nagata}(4)).
Hence $x$ is contained in only finitely many
quasi--$\star_{_{\!f}}$--maximal ideals of $D$.  Therefore  we
conclude that    $D$ is a $\star$--DD.
\end{proof}

 From the previous result, we deduce immediately:

\begin{cor} \label{cor:4.19} Let $D$ be an integral domain. The following are equivalent:
\begin{enumerate}
\item[(1)] $D$ is a Krull domain.  \item[(2)] $\Na(D,v) \ (= \Kr(D, v))$
is a PID.  \item[(3)] $\Na(D, v) \ (= \Kr(D, v))$ is a Dedekind
domain.  \hfill $\Box$
\end{enumerate}
\end{cor}

 Another consequence of Theorem \ref{prop:3.7} is the following:

\begin{cor} \label{cor:3.8} Let $D$ be an integral domain and $\star$ a semistar
operation on $D$. The following are equivalent:
\begin{enumerate}
\item[(1)] $D$ is a $\star$--DD.
\item[(2)] $D$ is  a    $\tilde{\star}$--DD.
\item[(3)] $D^{\tilde{\star}}$
is a $\dot{\tilde{\star}}$-DD.
 \item[(4)] $D^{\tilde{\star}}$ is a Krull domain and
$\dot{\tilde{\star}} = t_{D^{\tilde{\star}}}$.
\end{enumerate}
\end{cor}
\begin{proof}
 The equivalence $(1) \Leftrightarrow (2) \Leftrightarrow (3)$
follows from Theorem \ref{prop:3.7} and the fact that
$\Na(D,\star)=\Na(D,{\tilde{\star}})=\Na(D^{{\tilde{\star}}},\dot{\tilde{\star}})$
  (Proposition \ref{prop:nagata}(6)).    \\
 The equivalence $(3) \Leftrightarrow (4)$ follows from Proposition
\ref{prop:krull}, using the fact that $(1) \Leftrightarrow (3)$.
\end{proof}

\begin{rem} \label{rk:4.17}   From Corollary \ref{cor:3.8}
((1)$\Leftrightarrow$(4)), we have that if $D$ is a $\star$--DD then
$T:= D^{\tilde{\star}}$ is a Krull domain and $\tilde{\star}
=$ \d{$({t}_T)$}$_{_{\!D}}$ (where
$t_T$ is the $t$--operation of $T$). 
Note that it is not true in general that, if $T$ is a Krull overring of an
integral domain $D$, then $D$ is a
\d{$({t}_T)$}$_{_{\!D}}$--Dedekind   domain .  \\
For instance, let $K$ be a field  and 
$X$ an indeterminate over $K$.  Set $T := K[\![X]\!]$, $ M := XT$ and $D
:= K[\![X^{2}, X^{3} ]\!]  $.  It is easy to see that $D$ is a
one-dimensional local Noetherian integral domain with integral closure
equal to $T$ and maximal ideal equal to $N := M \cap D$ (with $NT =N$).
Therefore, in this case, ${t}_T = d_{T}$ is the identity (semi)star
operation on $T$ and so the semistar operation    \d{$({t}_T)$}$_{_{\!D}}$  on $D$
coincides with $\star_{\{T\}}$.  Clearly $\widetilde{\  \star_{\{T\}}} =
d_{D} \ (\lneq \star_{\{T\}})$, since $\mathcal{M}({\star_{\{T\}}})
=\{N\}$ and, obviously $D = D_{N}$ is not a DVR. Therefore $D$ is not a
 \d{$({t}_T)$}$_{_{\!D}}$--Dedekind domain.  \rm

 From the positive side, we have the following answer to the
question of when, given a Krull overring $T$ of an integral domain
$D$, $D$ is a \d{$({t}_T)$}$_{_{\!D}}$--DD:

\bf (\ref{rk:4.17}.1) \it Let $T$ be an overring of an integral
domain $D$. The following are equivalent:
\begin{enumerate}
\item[(1)] $D$ is a \d{$({t}_T)$}$_{_{\!D}}$--DD.
\item[(2)] $T$ is a
Krull domain  and, for each
$t_T$--maximal ideal $Q$ of $T$, $D_{Q \cap D}= T_Q $.
\end{enumerate} \rm

The previous characterization is a straightforward consequence of the
following ``restatement'' of the
equivalence given in Corollary \ref{cor:3.8}
((1)$\Leftrightarrow$(4)): 

\bf (\ref{rk:4.17}.2) \it If $D$ is an integral domain and $\star$ is
a semistar operation on $D$, then the following are equivalent:
\begin{enumerate}
\item[(1)] $D$ is a $\star$--DD.
\item[(2)] There exists an overring $T$ of $D$ such that $T$ is a
Krull domain,   $\star_{_{\!f}} = $   \d{$({t}_T)$}$_{_{\!D}}$   and, for each
$t_T$--maximal ideal $Q$ of $T$, $D_{Q \cap D}= T_Q $.
\end{enumerate} \rm

 To show the previous equivalence, note that in general the set of the
quasi--\d{$({t}_T)$}$_{_{\!D}}$--maximal ideals in $D$ coincide with
the set $\{Q\cap D \mid Q \mbox{ is a $t_{T}$--maximal ideal in
$T$ }\}$ \cite[Lemma 2.3(3)]{FL03}. Therefore the assumption that
$\star_{_{\!f}} = $ \d{$({t}_T)$}$_{_{\!D}}$  and, for each
$t_T$--maximal ideal $Q$ of $T$, $D_{Q \cap D}= T_Q $ implies that
$E^{\tilde{\star}} =(ET)^{t_{T}}$, for each $E \in
{\overline{\boldsymbol{F}}}(D)$, (in particular, $D^{\tilde{\star}}
=T$), and so $\dot{\tilde{\star}}^T = t_{T}$. Therefore
(\ref{rk:4.17}.2(2)) implies condition (4) of  Corollary
\ref{cor:3.8}.\\
Conversely, assume that  condition (4) of  Corollary
\ref{cor:3.8} holds and set $T :=D^{\tilde{\star}}$. Since
$\dot{\tilde{\star}} = t_{T}$ and $D$ is a
$\star$--DD (Corollary
\ref{cor:3.8} ((4)$\Rightarrow$(1))), then  $\tilde{\star}=
\star_{_{\!f}}$ (Proposition \ref{pr:1.5} and \ref{prop:3.2}), and so
$\dot{{\star}}_{_{\!f}}  = t_{T}$. Therefore
${\star}_{_{\!f}}  = $ \d{$({t}_T)$}$_{_{\!D}}$ (Proposition
\ref{pr:1.2}(2)). Moreover, by the previous considerations,
the set of the
quasi--${\star}_{_{\!f}} $--maximal ideals in $D$ coincide with
the set $\{Q\cap D \mid Q \mbox{ is a $t_{T}$--maximal ideal in
$T$ }\}$. Since $D$ is a $\star$--DD and hence, in particular, a
$\star$--ADD (and since $T$ is a Krull domain), then $D_{Q \cap D}$
is a DVR, which must coincide with its (DVR) overring
$T_Q $, for each $t_T$--maximal ideal $Q$ of $T$.

 It is possible to give another proof of (\ref{rk:4.17}.2) by using
Lemma \ref{lemma:noeth}(2)) and showing
the following   preliminary result of intrinsical interest
concerning the P$\star$MDs:

\bf (\ref{rk:4.17}.3) \it Let $D$ be an integral domain and $\star$ a
semistar operation on $D$.  Then, the following are equivalent:
\begin{enumerate}
\item[(1)] $D$ is a P$\star$MD.
\item[(2)] There exists an overring $T$ of $D$ such that $T$ is a
P$v_{T}$MD, $\star_{_{\!f}} =$   \d{$({t}_T)$}$_{_{\!D}}$  and, for each
$t_T$--maximal ideal $Q$ of $T$, $D_{Q \cap D}= T_Q $.
\end{enumerate}  \rm

 The proof is based on a variation of the techniques already discussed above
and the details are omitted. 

\end{rem}

\begin{ex} \label{ex:4.22} \sl  Let $D$ be a Mori domain, let $\Theta$ be the set of
all the maximal $t$--ideals of $D$ which are $t$--invertible and let
$\star_{\Theta}$ be the spectral semistar operation on $D$ associated
to $\Theta$ (Example \ref{ex:1.1} (3)). Assume that $\Theta \neq
\emptyset$ (i.e. that $D$ is a Mori non strongly Mori domain,
accordingly to the terminology introduced by Barucci and Gabelli \cite[page 105]{BG87}), then
$D$ is a $\star_{\Theta}$--DD. \rm

We apply the characterization given in (\ref{rk:4.17}.1) or in
Corollary \ref{cor:3.8} ((1)$\Leftrightarrow$(4)).  Note that
by \cite[Proposition 3.1 and Theorem 3.3 (a)]{BG87},
$D^{\star_{\Theta}}$
is a Krull domain such that the map $P \mapsto P^{\star_{\Theta}}$
defines a bijection between $\Theta$ and the set
$\mathcal{M}(t_{D^{\star_{\Theta}}})$ of all the
$t$--maximal ideals of  $D^{\star_{\Theta}}$ and $D_{P}=
(D^{\star_{\Theta}})_{P^{\star_{\Theta}}}$.  Therefore the (semi)star
operation $\dot{\star}_{\Theta}$ on $D^{\star_{\Theta}}$ coincides
with the $t$--operation, $t_{D^{\star_{\Theta}}}$, on
$D^{\star_{\Theta}}$.  Moreover, it is easy to see that, on $D$,  the semistar
operation
$(\mbox{\d{t}}_{D^{\star_{\Theta}}})_{_{\!D}}$ coincides with $\star_{\Theta}$.
\end{ex}

\vskip 12pt

Let $D$ be an integral domain and $\star$ a semistar operation on
$D$. We say that two nonzero ideals $A$ and $B$ are
\emph{$\star$--comaximal} if $(A+B)^\star=D^\star$. Note that, if
$\star$ is a semistar operation of finite type, then $A$ and $B$
are $\star$--comaximal if and only if $A$ and $B$ are not contained
in a common quasi--$\star$--maximal ideal.

\begin{lemma} \label{lemma:GP4.2}
Let $D$ be an integral domain and $\star$ a semistar operation on
$D$.  Let $A$ and $B$ be two  nonzero    $\star$--comaximal ideals of
$D$.  Then $(A \cap B) ^\star = (AB)^\star$.
\end{lemma}
\begin{proof}
In general $(A+B)(A \cap B) \subseteq AB$. Then, $((A+B)(A \cap
B))^\star \subseteq (AB)^\star \subseteq (A \cap B)^\star$. But
$((A+B)(A \cap B))^\star=((A+B)^\star(A \cap B))^\star=(D^\star (A
\cap B))^\star = (A \cap B)^\star$. Hence, $(A \cap B)^\star =
(AB)^\star$.
\end{proof}

\begin{cor} \label{cor:GP4.3}
Let $D$ be an integral domain and $\star$ a semistar operation of
finite type.  Let $n \geq 2$ and let $A_1, A_2, \ldots, A_n$ be  nonzero
ideals of $D$, such that $(A_i + A_j)^{\star} = D^{\star}$, for $i \neq
j$.  Then, $(A_1 \cap A_2 \cap \ldots \cap A_n)^{\star} = (A_1A_2
\cdot\ldots\cdot
A_n)^{\star}$.
\end{cor}
\begin{proof}
We prove it by induction  on $n \geq 2$,    using Lemma
\ref{lemma:GP4.2}  for the case $n =2$.     Set  $A :=    A_1 \cap
A_2 \cap \ldots \cap A_{n-1}$ and  $B:=   A_n$.  Then, $A$ and
$B$ are not contained in a common quasi--$\star$--maximal ideal,
otherwise, $A_n$ and $A_j$ (for some $1 \le j \le n-1$) would be contained
in a common quasi--$\star$--maximal ideal.  Hence $(A_1 \cap A_2 \cap
\ldots \cap A_{n-1} \cap A_n)^{\star} = (A \cap B)^{\star}
=(AB)^{\star}=(A^{\star}B)^{\star}=(A_1A_2 \cdot\ldots\cdot A_n)^{\star}$.
\end{proof}

\begin{thm} \label{thm:3.6} Let $D$ be an integral domain and $\star$ semistar
operation. The following are equivalent:
\begin{enumerate}
\item[(1)] $D$ is a $\star$-DD.
\item[(2)] For each nonzero ideal $I$ of $D$, there exists
a finite family of
quasi--$\star_{_{\!f}}$--prime    ideals $P_1, P_2,\ldots, P_n$  of
$D$, pairwise $\star_{_{\!f}}$-comaximals,  and a finite family of non
negative integers $e_1, e_2,\ldots, e_n$    such that
$I^{\tilde{\star}}= (P^{e_1}_1P^{e_2}_2\cdot\ldots\cdot P^{e_n}_n)^{\tilde{\star}}
$.
\end{enumerate}
 Moreover, if (2) holds and if $I^{\tilde{\star}} \neq D^{\tilde{\star}}$, then we can
assume that ${P_i}^{\tilde{\star}} \neq D^{\tilde{\star}}$, for each
$i=1, 2,\ldots, n$.  In this case, the integers $e_1, e_2,
\ldots, e_n$ are positive and the factorization is unique. 
\end{thm}
\begin{proof}
(1) $\Rightarrow$ (2).  Let $I$ be a nonzero ideal of $D$.   To avoid the
trivial case, we can assume that $I^{\tilde{\star}} \neq
D^{\tilde{\star}}$.     Let $P_1, P_2, \dots, P_n$ be the finite  (non
empty)    set of quasi--$\star_{_{\!f}}$--maximal ideals such that $I\subseteq
P_i$,  for $1 \leq i \leq n$    (Theorem \ref{thm:3.3}).  We have
$I^{\tilde{\star}} =\cap\{ID_P\mid\, P\in \mathcal{M}(\star_{_{\!f}})\}
= \cap_{i=1}^{i=n}(ID_{P_i} \cap D^{\tilde{\star}} )$.  Since $D_{P_i}$
is a DVR, then $ID_{P_i}=P_i^{e_i}D_{P_i}$, for some
integers    $e_i\ge 1$, $i=1, 2,\dots, n$.  Therefore, we have $ID_{P_i}\cap
D^{\tilde{\star}}=  P_i^{e_i}D_{P_i} \cap D^{\tilde{\star}} =
(P_i^{e_i})^{\tilde{\star}}$.  Hence
$I^{\tilde{\star}}=(P_1^{e_1})^{\tilde{\star}} \cap
(P_2^{e_2})^{\tilde{\star}} \cap \ldots \cap (P_n^{e_n})^{\tilde{\star}}
= (P^{e_1}_1 \cap P^{e_2}_2 \cap \ldots \cap P^{e_n}_n)^{\tilde{\star}}=
(P^{e_1}_1 P^{e_2}_2\cdot\ldots\cdot P^{e_n}_n)^{\tilde{\star}}$, by Corollary
\ref{cor:GP4.3}.\\
 For the  last statement,    let
$I^{\tilde{\star}}=(P^{e_1}_1P^{e_2}_2\cdot\ldots\cdot P^{e_n}_n)^{\tilde{\star}}$, if
${P_i}^{\tilde{\star}} = D^{\tilde{\star}}$, for some $i$, then obviously we can
cancel $P_i$ from the factorization of $I^{\tilde{\star}}$.   \\
 We
prove next  the uniqueness of the representation of $I^{\tilde{\star}}$.
  From  (Proposition \ref{prop:nagata}(4)), 
  we deduce that
  $I\!\Na(D,\star)=P^{e_1}_1 P^{e_2}_2\cdot\ldots\cdot P^{e_n}_n\!\Na(D,\star)=$
  $(P_1\!\Na(D,\star))^{e_1}$ $(P_2\!\Na(D,\star))^{e_2} \cdot\ldots\cdot
  (P_n\!\Na(D,\star))^{e_n}$ is the unique factorization into primes of the
  ideal $I\!\Na(D,\star)$ in the PID $\Na(D,\star)$ (Theorem
  \ref{prop:3.7}).  Since $P_i =P_i\!\Na(D,\star) \cap D$  (because each
  $P_i $ is a  quasi--$\widetilde{\star}$--maximal  ideal of $D$),    the
  factorization of $I^{\tilde{\star}}$ is unique.

  \noindent (2) $\Rightarrow$ (1)  Without loss of generality, we can
  assume that $D$ is not a field.     First, we prove that each
  localization to a quasi--$\star_{_{\!f}}$--maximal ideal of $D$ is a
   DVR. Let  $M    \in \mathcal{M}(\star_{_{\!f}})$ and let $J$
  be a  nonzero    proper ideal of  $D_M$.      Set $I:=    J \cap
  D \  ( \subseteq M )$.     Then,  it is easy to see that
  $I^{\tilde{\star}} \neq D^{\tilde{\star}}$ thus, by assumption,
  $I^{\tilde{\star}} = (P^{e_1}_1 P^{e_2}_2\cdot\ldots\cdot
  P^{e_n}_n)^{\tilde{\star}}$, for some  family of
  quasi--$\star_{_{\!f}}$--prime ideals $P_i$, with ${P_i}^{\tilde{\star}}
  \neq D^{\tilde{\star}}$ and for some family of integers $e_i \ge 1$,
  $i=1, 2, \dots, n$.
  It
  follows that  $J= ID_M    = I^{\tilde{\star}} D_M   =
  ({P_1}^{e_1} P^{e_2}_2\cdot\ldots\cdot {P_n}^{e_n})^{\tilde{\star}} D_M
    =(P^{e_1}_1 P^{e_2}_2\cdot\ldots\cdot  P^{e_n}_n) D_M$   (since
  $\tilde{\star}$ is a spectral semistar operation  defined by the set
  $\mathcal{M}(\star_{_{\!f}})$). Hence    $J$ is  a finite product of
  primes of  $D_{M}$.
  Therefore
   $D_{M}$    is a local Dedekind domain   \cite[Theorem
  37.8 ((1)$\Leftrightarrow$(3))]{Gil92},    that is,  $D_{M}$ is    a DVR.\\
Now we show that each
quasi--$\tilde{\star}$--prime ideal  of $D$    is
quasi--$\tilde{\star}$--invertible.  Let  $Q$    be a
quasi--$\tilde{\star}$--prime  of $D$    and let $0 \neq  x \in Q$.
Then,  by assumption, $(xD)^{\tilde{\star}}   = (P^{e_1}_1
P^{e_2}_2\cdot\ldots\cdot P^{e_n}_n)^{\tilde{\star}}$,   with $P_1, P_{2}, \ldots,
P_n$ nonzero prime ideals of $D$  and $e_i \ge 1$, $i=1, 2, \dots,
n$.     Since  $xD$     is obviously invertible (and thus, clearly,
quasi--$\tilde{\star}$--invertible), then each $P_i$ is
quasi--$\tilde{\star}$--invertible \cite[Lemma 2.10]{FP}.  Moreover,
since  $Q$    is a quasi--$\tilde{\star}$--ideal  of $D$,   then
$P^{e_1}_1 P^{e_2}_2\cdot\ldots\cdot P^{e_n}_n \subseteq (P^{e_1}_1 P^{e_2}_2\cdot\ldots\cdot
P^{e_n}_n)^{\tilde{\star}} \cap D \subseteq  Q$.  Therefore, $P_j
\subseteq  Q$    for some  $j$, with $1\leq j\leq n$, and
since  $D_Q$   is a DVR, we have  $Q= P_j$.     Hence
 $Q$     is a quasi--${\tilde{\star}}$--invertible  ideal of $D$.
Therefore,  by Theorem \ref{thm:GP3.4}, we conclude that    $D$ is
$\tilde{\star}$--Dedekind.
\end{proof}

 \begin{rem}
 It is clear that, if $D$ is a $\star$--DD then, for each nonzero
 ideal $I$ of $D$, such that $I^{\star_{_{\!f}}} \neq
 D^{\star_{_{\!f}}}$, we have a unique  factorization
 $I^{\star_{_{\!f}}}=(P^{e_1}_1 P^{e_2}_2\cdot\ldots\cdot
 P^{e_n}_n)^{\star_{_{\!f}}}$,  for some family of
 quasi--$\star_{_{\!f}}$--prime ideals $P_i$, with $P_i^{\star_{_{\!f}}}
 \neq D^{\star_{_{\!f}}}$, and for some family of positive integers $e_i $,
 $i=1, 2, \dots, n$,    since $\tilde{\star} = \star_{_{\!f}}$
   (Proposition \ref{pr:1.5}).
 The converse is
 not true.  For instance, take $D$, $T$ and $\star$ as in Remark
 \ref{rem:2.1}(2).    For  each nonzero proper ideal    $I$ of $D$, we
 have $I^{\star_{_{\!f}}}=IT=M^e=(M^e)^{\star_{_{\!f}}}$,  for some
 positive integer $e$,    since $T$ is a     DVR.   Note
 that this representation is unique, since $D$ is local with maximal
 ideal $M$ and $\di(D) =1$.
  But we
 have already observed that $D$ is not a $\star$--DD.
\end{rem}

 Next result generalizes to the semistar setting \cite[Theorem
38.5 ((1)$\Leftrightarrow$(3))]{Gil92}.

\begin{thm} \label{thm:3.9}
Let $D$ be an integral domain  which is not a field    and $\star$ a
semistar operation on $D$.  The following are equivalent:
\begin{enumerate}
\item[(1)] $D$ is a $\star$--DD.
\item[(2)] For each nonzero ideal $I$ and for each $a\in
I$, $a\neq 0$, there exists  $b\in I^{\widetilde{\star}}$  such that
$I^{\tilde{\star}}=((a, b)D)^{\tilde{\star}}$.
\end{enumerate}
\end{thm}
\begin{proof}
(1) $ \Rightarrow $ (2).    We start by proving the following:

\bf Claim. \sl  If $D$ is a $\star$--DD, then the map $M \mapsto
M^{\widetilde{\star}}$ establishes a bijection between the set
$\mathcal{M}(\star_{_{\!f}})$ \ ($=  \mathcal{M}(\widetilde{\star})$
by Proposition \ref{prop:nagata} (5))
of the
quasi--$\star_{_{\!f}}$--maximal ideals of $D$ and the set
$\mathcal{M}(t_{D^{\widetilde{\star}}})$ of the
$t_{D^{\widetilde{\star}}}$--maximal ideals of (the Krull domain)
$D^{\widetilde{\star}}$. \rm

For each $M \in \mathcal{M}(\star_{_{\!f}})$, since
${D^{\widetilde{\star}}} \subseteq D_{M}$,  it is easy to see that
$M^{\widetilde{\star}} = MD_{M}\cap  {D^{\widetilde{\star}}}$.
Therefore,
$M^{\widetilde{\star}}$ is a $\dot{\widetilde{\star}}$--prime ideal
of  $D^{\widetilde{\star}}$ and  $M^{\widetilde{\star}} \cap D = M$.
Furthermore, by Corollary \ref{cor:3.8}, we know
that  $D^{\widetilde{\star}}$ is a Krull domain and
$\dot{\widetilde{\star}} = t_{D^{\widetilde{\star}}}$.  On the other
hand, for
each $\dot{\widetilde{\star}}$--prime ideal $N$ of
$D^{\widetilde{\star}}$, we know that $N\cap D$ is a
quasi--${\widetilde{\star}}$--prime of $D$ \cite[Lemma 2.3 (4)]{FL03}).
Since $D$ is a
$\star$--DD (or, equivalently, a ${\widetilde{\star}}$--DD),  we have that each
quasi--${\widetilde{\star}}$--prime is a
quasi--${\widetilde{\star}}$--maximal  (Proposition \ref{prop:2.2} (2)),
thus we easily conclude.

Let $a\in I$, $a\neq 0$, and $\{M_1, M_2,
\ldots,  M_n  \}$ the (finite) set of
quasi--$\star_{_{\!f}}$--maximal ideals such that $a\in M_i$.  Since
$D_{M_i}$ is a DVR, then $ID_{M_i}=x_iD_{M_i}$, for some $x_i
\in I$,  for each $i=1, 2, \ldots,  n$.    We use the fact that
$D^{\widetilde{\star}}$ is a Krull domain and, by the Claim, that
$\{D^{\widetilde{\star}}_{M^{\widetilde{\star}}}= D_{M} \mid M \in
\mathcal{M}(\star_{_{\!f}})\}$ is the defining family of the rank-one
discrete valuation overrings of $D^{\widetilde{\star}}$, in order to
apply the
approximation theorem to $D^{\widetilde{\star}}$.   Let $v_1, v_2, \ldots, v_n$ be the
valuations associated  respectively    to $D_{M_1}, D_{M_2}, \ldots,
 D_{M_n}$ and let $v_{M'}$ be the valuation associated to $D_{M'}=
 D^{\widetilde{\star}}_{{M'}^{\widetilde{\star}}}$,  for $M' \in
\mathcal{M'}:= \mathcal{M}(\star_{_{\!f}}) \setminus \{M_1, M_2, \ldots, M_n\}$.  Set $k_1:=v_1(x_1), k_2:=v_2(x_2), \ldots,
 k_n:= v_n(x_n)$.       Then there exists $b\in K$ such that $v_i(b)=k_i$,
for each $i=1, 2, \ldots,  n$,   and $v_{M'}(b) \geq 0$,  for each $M'
\in \mathcal{M'}$ \cite[Theorem 44.1]{Gil92}.   
We have $I^{\tilde{\star}}=((a, b)D)^{\tilde{\star}}$.  Indeed, let
$M\in \mathcal{M}(\star_{_{\!f}})$.  If $M=M_i$, for some $i$, then
$ID_{M}=ID_{M_i}=x_iD_{M_i}=bD_{M_i}=(a, b)D_{M_i}$.  If $M\neq M_i$ for
each $i$, then $ID_{M}=D_M=(a, b)D_M$.\\
(2) $ \Rightarrow $ (1).  Let $M\in \mathcal{M}(\star_{_{\!f}})$ and
$J$     a nonzero ideal of $D_M$.  Let $a\in  J$,   $a\neq 0$,    there exists $s\in D$,
$s\notin M$, such that $sa\in  I:=J\cap D$.    Then, by assumption, there
exists  $b \in  I^{\widetilde{\star}}$     such that  $I^{\tilde{\star}}    = ((sa,
b)D)^{\tilde{\star}}$.   Therefore,    we have
 $J=ID_M=I^{\tilde{\star}}D_M    =((sa, b)D)^{\tilde{\star}}D_M = (sa,
b)D_M=(a, b)D_M$.  By \cite[Theorem 38.5]{Gil92}, $D_M$ is a Dedekind
domain, and hence a DVR. Thus, $D$ is a $\star$--ADD,  hence, in
particular, is is a P${\tilde{\star}}$MD
(Corollary \ref{cor:2.6} 
and Proposition \ref{prop:2.2}(2)).    In addition,  from the
assumption  and from \cite[Lemma 2.3]{FP}, we deduce that  
$D$ is $\tilde{\star}$--Noetherian (Lemma \ref{lemma:3.1}), hence
$D$ is a $\star$--DD  (Corollary \ref{cor:tilde} and Proposition
\ref{prop:3.2}(1)).
\end{proof}

\begin{rem}
 Note that,   if $D$ is a $\star$--DD  (and hence
 $\widetilde{\star} ={\star_{_{\!f}}}$),  then $D$ satisfies also
  a
 statement concerning ${\star_{_{\!f}}}$, analogous to the statement (2)
 in Theorem \ref{thm:3.9}: 
\begin{enumerate}
 \item[($2_{_{\!f}}$)]  for each nonzero ideal $I$ of $D$ and for each
$0 \neq a \in I$, there exists $b \in  I^{\star_{_{\!f}}}$  such that
$((a,b)D)^{\star_{_{\!f}}}=I^{\star_{_{\!f}}}$.
\end{enumerate}
But  ($2_{_{\!f}}$)    does not imply that    $D$  is a
$\star$--DD.  For instance,    let $D, T$ and $\star$  be    as in
Remark \ref{rem:2.1}.  Obviously, for each nonzero  proper     ideal
$I$ of $D$ and    for each nonzero $a \in I \subseteq D$ we have
$I^{\star_{_{\!f}}}=IT=  X^nT = (a, X^n)T=((a,
X^n)D)^{\star_{_{\!f}}}$, for some $n \geq 1$,  (where $X^n \in
I^{\star_{_{\!f}}}\cap D$),    but $D$ is not a
$\star$--DD.
\end{rem}


\end{document}